\documentclass[11pt]{amsart}
\usepackage[margin=1in]{geometry}
\usepackage{amsmath,amsthm,mathtools,enumitem,amssymb}
\usepackage{newtxtext}
\usepackage[alphabetic,nobysame]{amsrefs}
\renewcommand{\MR}[1]{} \renewcommand{\PrintDOI}[1]{}
\usepackage[mathscr]{euscript}
\usepackage[colorlinks=true,linkcolor=blue,citecolor=blue]{hyperref}
\hypersetup{
  pdftitle={Classification of twisted Yangians as coideal subalgebras of Yangians},
  pdfkeywords={Yangian, twisted Yangian, coideal quantization,
    quantum symmetric pair, minimalistic presentation}
}

\numberwithin{equation}{section}
\newtheorem{theorem}{Theorem}[section]
\newtheorem{proposition}[theorem]{Proposition}
\newtheorem{lemma}[theorem]{Lemma}
\newtheorem{corollary}[theorem]{Corollary}
\newtheorem{maintheorem}{Theorem}

\theoremstyle{definition}
\newtheorem{definition}[theorem]{Definition}
\theoremstyle{remark}
\newtheorem{remark}[theorem]{Remark}
\newtheorem{example}[theorem]{Example}

\newcommand{\C}{\mathbb C}
\newcommand{\B}{{\mathbf B}}
\newcommand{\K}{{\mathscr K}}
\newcommand{\I}{{\mathbb I}}
\newcommand{\hI}{{\widehat{\mathbb I}}}
\newcommand{\rU}{{\mathrm U}}
\newcommand{\rS}{{\mathrm S}}
\newcommand{\g}{\mathfrak g}
\newcommand{\kalg}{\mathfrak k}
\newcommand{\palg}{\mathfrak p}
\newcommand{\Y}{\mathrm Y_{\hbar}(\g)}
\newcommand{\Ltw}{\mathfrak g[u]^{\check\omega}}
\newcommand{\cl}{\mathrm{cl}}
\newcommand{\cli}{\mathrm{cl}_{\imath}}
\newcommand{\inc}{\mathrm{inc}}
\newcommand{\ad}{\operatorname{ad}}
\newcommand{\Hom}{\operatorname{Hom}}
\newcommand{\End}{\operatorname{End}}
\newcommand{\pr}{\operatorname{pr}}
\newcommand{\sfb}{\mathsf b}
\newcommand{\arxiv}[1]{\href{http://arxiv.org/abs/#1}{arXiv:\nolinkurl{#1}}}

\title{Classification of twisted Yangians as coideal subalgebras of Yangians}
\author{Kang Lu}
\address{
Shenzhen International Center for Mathematics and Department of Mathematics, Southern University of Science and Technology, Shenzhen, China
}
\email{kang.lu.math@gmail.com,luk@sustech.edu.cn}
\begin{document}

\begin{abstract}
Let $\mathfrak g$ be a complex simple Lie algebra with an involution $\omega$. We classify the homogeneous right coideal subalgebras of the Yangian of $\mathfrak g$ whose classical limit is the enveloping algebra of the twisted current algebra, embedded in the standard way, without assuming any formula for the coaction. For nontrivial $\omega$, they form a family parametrized by a spectral parameter and an element of the center of the fixed-point subalgebra; for $\omega=\mathrm{id}$, only a spectral parameter occurs. All members of the family are isomorphic as graded algebras, through explicit isomorphisms. In particular, the twisted Yangians with a central correction introduced in \cite{HL26} are coideal quantizations, and the $\imath$Yangian of quasi-split type $\mathsf A_{2n}$ is isomorphic to the traditional twisted Yangian in \cite{Mac02,BR17}. The PBW property is proved uniformly from the classical affine Serre relations, without quantum defining relations. We also recover finite presentations in split and quasi-split types and verify the technical hypotheses of the quantum duality principle used in \cite{LWW25}.
\end{abstract}


\maketitle


\section{Introduction}\label{sec:introduction}

\subsection{Background}

Yangians are the Hopf algebra quantizations of the rational Lie
bialgebra structure on the current algebra $\g[u]$. Drinfeld
introduced them in connection with the quantum Yang--Baxter equation
\cite{Dri85,Dr87b} and gave their $\mathcal J$- and current
presentations \cite{Dri85,Dri87}. The current presentation underlies
the classification of finite-dimensional irreducible representations
\cite{CP91,CP94}, the $R$-matrix presentation is the natural framework
for the classical Lie algebras \cite{Mo07}, and the equivalence of the
three presentations is treated in \cite{BK05,JLM18,GRW19,Wen18}.

Twisted Yangians are coideal subalgebras of Yangians associated with
symmetric pairs. Their $R$-matrix realizations originate in the
reflection equation \cite{Che84,Skl88}. They were introduced by
Olshanski for types $\mathsf{AI}$ and $\mathsf{AII}$
\cite{Ols92,MNO96} and constructed for type $\mathsf{AIII}$
\cite{MR02} and for types $\mathsf B$, $\mathsf C$ and $\mathsf D$
\cite{GR16}, and their finite-dimensional irreducible representations
were classified in these settings \cite{Mol98,MR02,GRW17rep,GRW19rep}.
For general symmetric pairs, MacKay \cite{Mac02} constructed coideal
generators in Drinfeld's $\mathcal J$-presentation. Belliard and
Crampé \cite{BC12} obtained the twisted Yangians of $\mathfrak{sl}_2$
by deforming Lie bi-ideal structures. The $\mathcal J$-approach was
developed systematically by Belliard and Regelskis \cite{BR17}. In
particular, they computed the Lie coideal structure of the twisted current algebra
\cite[Props.~3.2 and~4.5]{BR17}, the coaction on Casimir-corrected
degree-one generators \cite[\S7.3 and~\S7.5]{BR17}, and, for the
identity involution, the degree-two generators \cite[\S7.4]{BR17}.
Current presentations were obtained for type $\mathsf{AI}$
\cite{LWZ25GD}, for all split types by degeneration of affine
$\imath$quantum groups \cite{LWZ25affine}, and
for quasi-split types \cite{LZ24}; see \cite{regelskis2026twisted} for
their comparison with $R$-matrix presentations in types
$\mathsf{BI}$, $\mathsf{CI}$ and $\mathsf{DI}$. Finite presentations
with generators in degrees zero and one were given for Yangians in
\cite{Lev93gen,GNW18} and for twisted Yangians of split and
nontrivial quasi-split types in \cite{Lu26min,HL26}.

\subsection{Motivation}

This paper is motivated by three developments.

One of the main motivations for this work comes from the twisted
Yangians with a central correction introduced in \cite{HL26}. In the
$\mathcal J$-presentation, a twisted Yangian is traditionally
generated by $\kalg$ and the elements
$\mathcal J(y)-\tfrac14\hbar[y,C_{\kalg}]$ for $y\in\palg$
\cite{Mac02,BR17}, where $\mathcal J$ denotes Drinfeld's degree-one
generators and $C_{\kalg}$ is the quadratic Casimir element of
$\kalg$. In \cite[Def.~2.8]{HL26}, these generators are replaced, in
our normalization, by
\[
 \mathcal J(y)-\tfrac14\hbar[y,C_{\kalg}]-c\hbar[y,\mathcal K]
 \qquad(y\in\palg),
\]
where $c\in\C$ and $\mathcal K$ lies in the center $\mathfrak z(\kalg)$.
This modification is motivated by an example in quasi-split type
$\mathsf A_{2n}$: the $\imath$Yangian in Drinfeld type presentation
corresponds to $c=\tfrac18$ and a specific nonzero $\mathcal K$,
rather than to the traditional twisted Yangian
\cite[Thm.~3.4(3)]{HL26}. It is then natural to ask which subalgebras
of this kind are coideal quantizations of the twisted current
algebra, whether there are others, and how the resulting algebras are
related to each other.

A second motivation comes from the uniqueness of Yangians and the
notion of coideal quantizations. Gautam, Wendlandt and Xu
\cite{GWX25} proved that the rational Lie bialgebra $\g[u]$ admits at
most one homogeneous quantization; their proof constructs normalized
degree-one generators in an arbitrary Hopf algebra quantization. On
the coideal side, a coideal quantization of the twisted current
algebra is defined in \cite[Def.~3.15]{LWW25} by prescribing only its
coideal property and its embedded classical limit. In this framework,
it is shown in \cite[Thm.~3.19]{LWW25}, via the quantum duality
principle \cite{Dr87b,G02,CG06}, that the quantum duality subalgebra
of a twisted Yangian quantizes the Poisson symmetric space
$G_1[\![u^{-1}]\!]/G_1[\![u^{-1}]\!]^{\omega}$, where
$G_1[\![u^{-1}]\!]$ is the first congruence subgroup of the loop
group. This is the unshifted counterpart of the quantization of
involutive fixed-point loci by shifted $\imath$Yangians. The argument
there uses the technical conditions \cite[(3.5) and (3.13)]{LWW25} on
lifts of bases. It is therefore natural to determine all coideal
quantizations in this sense inside a fixed Yangian, following the
strategy of normalized lifts in \cite{GWX25}.

A third motivation comes from the $\mathcal J$-approach of
\cite{BR17}. Existence and uniqueness up to isomorphism of
homogeneous quantizations of the twisted current algebra are stated
in \cite[Thms.~5.5 and~5.6]{BR17}, together with explicit generators
and quantum defining relations. The notion of quantization used there
is related to ours but is stricter: it imposes constraints on the
coproduct modulo $\hbar^2$, prescribing the first-order coaction by
the Lie coideal structure \cite[Def.~5.3]{BR17}; see
\cite[Rem.~3.17]{LWW25}. The flatness of these quantizations, via the
embedding into the Yangian, and their PBW theorem are discussed in
\cite[Rem.~5.8]{BR17}. In the present paper, we impose no condition on
the coaction beyond the coideal property. More precisely, the
first-order coaction is not part of the input; it is derived from the
coideal property and the classical limit (Lemma~\ref{lem:lift}). We
also do not start from chosen generators or relations. Both points
are essential if one wants to decide whether a given subalgebra, such
as those in \cite[Def.~2.8]{HL26}, is a coideal quantization, and to
determine all of them.

\subsection{Main results}\label{ssec:main-results}

The purpose of the present paper is to classify the homogeneous
coideal quantizations inside a fixed Yangian, for every involution of
every complex simple Lie algebra; see
Theorems~\ref{thm:introduction-main}
and~\ref{thm:introduction-identity} below. We show that the coideal
property alone determines the generators up to a spectral parameter
and, for Hermitian symmetric pairs, a central parameter. Every value
of the parameters occurs, and the resulting algebras are related by
explicit isomorphisms. In particular, we obtain uniqueness up to
graded algebra isomorphism, which recovers the uniqueness statements
of \cite[Thms.~5.5 and~5.6]{BR17} without prescribing the coaction.
The PBW property is proved by a saturation criterion
(Proposition~\ref{prop:saturation-criterion}), which uses only the
classical affine Serre relations, the coproduct and degree descent,
and no quantum defining relations, so that it applies to an arbitrary
coideal subalgebra.

Let $\g$ be a finite-dimensional complex simple Lie algebra, let
$\omega$ be an involution of $\g$, and write $\g=\kalg\oplus\palg$
for its decomposition into $(+1)$- and $(-1)$-eigenspaces. We fix the
invariant bilinear form $(\cdot,\cdot)$ on $\g$ normalized as in
\S\ref{ssec:symmetric-pairs}; Casimir elements and dual bases are
taken with respect to it. The automorphism
$\check\omega(xu^r)=(-1)^r\omega(x)u^r$ of $\g[u]$ has
the fixed-point subalgebra
\[
 \Ltw=\kalg[u^2]\oplus u\palg[u^2],
\]
the \emph{twisted current algebra}. We work in the homogeneous
Yangian $\Y$, with $\deg\hbar=\deg u=1$ and classical limit
$\Y/\hbar\Y\cong\rU(\g[u])$; see \S\ref{ssec:yangian}. A
\emph{homogeneous coideal quantization} of $\Ltw$, or a
\emph{twisted Yangian} associated with $(\g,\kalg)$, is a graded
unital subalgebra $B\subset\Y$, free as a graded $\C[\hbar]$-module,
such that
\[
 \Delta(B)\subset B\otimes\Y,
 \qquad
 B/\hbar B\cong\rU(\Ltw)\subset\rU(\g[u])
       \cong\Y/\hbar\Y,
\]
where the second identification is induced by the inclusion
$B\subset\Y$; see Definition~\ref{def:polynomial-coq}. This is the
right coideal version of \cite[Def.~3.15]{LWW25}; left coideals are
treated in \S\ref{ssec:left-coideal}. Note that neither a choice of
generators nor a formula for the coaction is part of the data. By the
Rees construction, all our results also apply to filtered right
coideal subalgebras of the ordinary Yangian $\mathrm Y(\g)$; see
Remark~\ref{rem:filtered}.

Suppose first that $\omega\ne\mathrm{id}$. For $a\in\C$ and
$\K\in\mathfrak z(\kalg)$, put
\begin{equation}\label{eq:introduction-generator}
 b_{a,\K}(y):=\mathcal J(y)-\tfrac14\hbar[y,C_{\kalg}]
       +\hbar\bigl(ay+[\K,y]\bigr),\qquad y\in\palg,
\end{equation}
and let $\B_{a,\K}$ be the $\C[\hbar]$-subalgebra of $\Y$ generated
by $\kalg$ and $b_{a,\K}(\palg)$. We call $a$ the \emph{spectral
parameter} and $\K$ the \emph{central parameter}. The spectral
parameter is shifted by the spectral translations
$\sigma_c:\mathcal J(x)\mapsto\mathcal J(x)+c\hbar x$ of $\Y$. The
central parameter can be nonzero only for Hermitian symmetric pairs,
since by Lemma~\ref{lem:modules}, $\mathfrak z(\kalg)$ is
one-dimensional if the $\kalg$-module $\palg$ is reducible, and zero
otherwise. In the Hermitian case, let $Z\in\mathfrak z(\kalg)$ act by
$\pm1$ on the two irreducible summands of $\palg$. Then $\K=cZ$
corrects $\mathcal J(y)-\tfrac14\hbar[y,C_{\kalg}]$ by $(a+c)\hbar y$
and $(a-c)\hbar y$ on these two summands, respectively.

\begin{maintheorem}[Nontrivial involutions]
\label{thm:introduction-main}
For every nontrivial complex symmetric pair $(\g,\kalg)$,
the following statements hold.
\begin{enumerate}[label=\textup{(\roman*)}]
\item For every $a\in\C$ and $\K\in\mathfrak z(\kalg)$,
$\B_{a,\K}$ is a homogeneous coideal quantization of $\Ltw$.
Ordered monomials in homogeneous lifts of an ordered
homogeneous basis of $\Ltw$ form a $\C[\hbar]$-basis
of $\B_{a,\K}$.
\item Every homogeneous coideal quantization $B\subset\Y$
of $\Ltw$ equals $\B_{a,\K}$ for a unique $a\in\C$ and a unique
$\K\in\mathfrak z(\kalg)$.
\item For all $a_0,a_1\in\C$ and $\K_0,\K_1\in\mathfrak z(\kalg)$, there
is a graded $\C[\hbar]$-algebra isomorphism
\[
 \Phi:\B_{a_0,\K_0}\xrightarrow{\ \sim\ }\B_{a_1,\K_1},\qquad
 x\longmapsto x+(\K_1-\K_0,x)1,\quad
 b_{a_0,\K_0}(y)\longmapsto b_{a_1,\K_1}(y),
\]
for $x\in\kalg$ and $y\in\palg$, such that
$\Delta\circ\Phi=(\Phi\otimes\sigma_{a_1-a_0})\circ\Delta$.
\end{enumerate}
\end{maintheorem}

Parts (i), (ii) and (iii) of Theorem~\ref{thm:introduction-main} are
proved in Theorems~\ref{thm:coideal-pbw},
\ref{thm:embedded-classification} and~\ref{thm:main-uniqueness},
respectively. The proof of (ii) shows that every homogeneous coideal
quantization has a unique normalized degree-one lift and its first-order coaction is fixed by the Lie coideal structure.
Moreover, as the embedded quantization varies, its normalized lift is
parametrized by $(a,\K)$. The isomorphism $\Phi$ in
(iii) is a composite of two kinds of isomorphisms. The first one is
the spectral translation $\sigma_{a_1-a_0}$, which is the restriction
of a Hopf algebra automorphism of $\Y$ and changes only $a$. The
second one is the central translation
$b\mapsto\chi(b_{(1)})\,b_{(2)}$, where $\chi$ is the character of
$\B_{a_0,\K_0}$ given by $(\K_1-\K_0,\cdot)$ on $\kalg$ and by zero on
$b_{a_0,\K_0}(\palg)$; it changes only $\K$. When $\K_0\ne\K_1$, the central translation changes the counit
on $\kalg$, therefore it is not the restriction of a Hopf algebra
automorphism of $\Y$. If $\K_0\ne\K_1$, then the
reduction of $\Phi$ modulo $\hbar$ is not the identity of $\rU(\Ltw)$,
but the automorphism induced by the character $(\K_1-\K_0,\cdot)$ of
$\kalg$. Therefore, the embedded twisted Yangians form a family
parametrized by $\C\times\mathfrak z(\kalg)$ (that is, by $\C^2$ for
Hermitian pairs and by $\C$ otherwise), although they are all
isomorphic as graded algebras.

As an application, the twisted Yangians of \cite[Def.~2.8]{HL26} are
the algebras $\B_{0,c\mathcal K}$; in particular, they are all coideal
quantizations of $\Ltw$. By \cite[Thm.~3.4(3)]{HL26}, the
$\imath$Yangian of quasi-split type $\mathsf A_{2n}$ corresponds to a
nonzero central parameter. Hence it is isomorphic as a graded algebra
to the traditional twisted Yangian $\B_{0,0}$; see
Corollary~\ref{cor:HL26}.

We mention that the central parameter has appeared before in special
cases. In type $\mathsf{AIII}$, it corresponds to the free parameter
of the $K$-matrix of the reflection algebras, and different values
give isomorphic algebras \cite[Def.~3.8, Prop.~3.4 and
Cor.~3.1]{CGM14}. For general Hermitian pairs, characters of the above
form appear as the counit in \cite[Thm.~5.5]{BR17}. By
Theorem~\ref{thm:introduction-main}, the spectral and central
parameters together account for all embedded twisted Yangians.

Parts (ii) and (iii) of Theorem~\ref{thm:introduction-main} have
counterparts for quantum symmetric pairs. The coideal subalgebras of
$U_q(\g)$ defining quantum symmetric pairs depend on parameters, and
Letzter characterized them as the coideal subalgebras that specialize
to $\rU(\kalg)$ in an appropriate maximal sense \cite{Let02}; see also
\cite{Kolb14,Molev2003coideal,DCNTY23}. Part (ii) can be viewed as a
Yangian analogue of this characterization. Moreover, $\imath$quantum
groups with different parameters are isomorphic
\cite[Lem.~2.5.1 and Rem.~2.5.3]{Wat21}. The isomorphism there is a
composite of the restriction of a Hopf algebra automorphism of the
quantum group, similar to our spectral translations, and a rescaling
of the Cartan generators of the coideal subalgebra, which does not
come from the quantum group, similar to our central translations.

For $\omega=\mathrm{id}$, the twisted current algebra is $\g[u^2]$,
which is generated in degrees zero and two. Let $\{x_r\}$ and
$\{x^r\}$ be dual bases of $\g$, let $C_\g=\sum_rx_rx^r$, and let
$c_\g\ne0$ be the eigenvalue of $C_\g$ on the adjoint representation.
The $\g$-module map
\[
 \mathcal G(x)=\frac1{c_{\g}}
       \sum_r[\mathcal J(x_r),\mathcal J([x^r,x])]
       -\tfrac14\hbar[\mathcal J(x),C_{\g}]
 \qquad(x\in\g)
\]
lifts $x\mapsto xu^2$ to $\Y_2$. Set
\[
 \mathcal G_a(x):=\mathcal G(x)+2a\hbar\mathcal J(x)
                         +a^2\hbar^2x,\qquad
 \B^{\mathrm{ev}}_a
 =\C[\hbar]\langle\g,\mathcal G_a(\g)\rangle\qquad(a\in\C).
\]

\begin{maintheorem}[The identity involution]
\label{thm:introduction-identity}
For every $a\in\C$, $\B^{\mathrm{ev}}_a$ is a homogeneous
coideal quantization of $\g[u^2]$, with a PBW basis obtained
by lifting an ordered homogeneous basis of $\g[u^2]$.
Every homogeneous coideal quantization of $\g[u^2]$ in
$\Y$ equals $\B^{\mathrm{ev}}_a$ for a unique $a\in\C$.
For $a_0,a_1\in\C$, the spectral translation $\sigma_{a_1-a_0}$ of
$\Y$ restricts to a graded algebra isomorphism
\[
 \sigma_{a_1-a_0}:\B^{\mathrm{ev}}_{a_0}
       \xrightarrow{\sim}\B^{\mathrm{ev}}_{a_1},\qquad x\mapsto x,\qquad  \mathcal G_{a_0}(x)\mapsto \mathcal G_{a_1}(x)\qquad (x\in\g).
\]
\end{maintheorem}

Theorem~\ref{thm:introduction-identity} is proved in
Theorems~\ref{thm:identity-coideal-pbw}
and~\ref{thm:identity-uniqueness}. Note that the degree-two lift
itself is not unique, but this ambiguity does not affect the
generated subalgebra. Our proofs of
Theorems~\ref{thm:introduction-main}
and~\ref{thm:introduction-identity} are uniform and do not use any
case-by-case analysis of symmetric pairs. Apart from general facts on
reductive Lie algebras, we only use Kac's realization of twisted loop
algebras, through the affine data in
Definition~\ref{def:affine-datum}.

\begin{remark}[Comparing constructions]\label{rem:compare}
By Theorems~\ref{thm:introduction-main}
and~\ref{thm:introduction-identity}, in order to compare different
constructions of twisted Yangians, it suffices to check that the
construction yields a graded right coideal subalgebra of $\Y$ and
that its natural classical limit is $\rU(\Ltw)\subset\rU(\g[u])$. Then
the parameters can be read off from the normalized lifts, and any two
such constructions give isomorphic graded algebras. This applies to
$R$-matrix, $\mathcal J$- and current realizations, as long as they
are realized inside the same Yangian of a simple Lie algebra.
\end{remark}

For split types other than $\mathsf A_1$ and nontrivial quasi-split
types other than $\mathsf A_2$, Theorem~\ref{thm:main-presentation}
gives a finite presentation with generators $\kalg$ and a degree-one
copy of $\palg$, whose only relations besides the Lie algebra and
module relations form a single family of Cartan commutators. This
recovers the minimalistic presentations of \cite{Lu26min,HL26} by a
direct argument in the $\mathcal J$-presentation. More precisely, the
relations are derived from the embedding, the completeness follows
from the presentation of the twisted current algebra, and all central
parameters are reached by central translations. Compared with the
$\mathcal J$-presentation of \cite{BR17}, only Cartan commutators are
needed; for symmetric pairs that are not quasi-split, we refer to
\cite{BR17} for the relations.

Finally, in Appendix~\ref{sec:quantum-duality-assumption} we show that
$B'=B\cap\Y'$ for the quantum duality subalgebras $\Y'\subset\Y$ and
$B'\subset B$, and that the technical conditions
\cite[(3.5) and (3.13)]{LWW25} hold for the Yangian and for every
twisted Yangian, for every involution and all parameters; see
Theorem~\ref{thm:duality-technical-assumption}. The statement for
Yangians is also recorded in \cite[Rem.~3.2]{LWW25}.

\subsection{Organization}

The paper is organized as follows. In Section~\ref{sec:conventions},
we fix the conventions and collect two tools used throughout the
paper, namely a relative primitive-element lemma and a graded lifting
lemma. In Section~\ref{sec:criterion}, we prove a saturation
criterion for subalgebras generated by $\rU(\kalg)$ and a lift of the
lowest positive-degree piece of $\Ltw$. The idea is that a first
failure of $\hbar$-saturation is detected by an affine Serre relation
whose weight and involution eigenvalue are excluded from the relevant
coefficient space, and then degree descent leads to a contradiction.
Combined with graded lifting, this yields PBW bases and a
presentation by lifted Serre relations. In Section~\ref{sec:rigidity}, we classify the
embeddings for $\omega\ne\mathrm{id}$. We normalize the coaction of a
degree-one lift by relative primitives, and the two parameters come
from the identity
$\Hom_{\kalg}(\palg,\g)=\C\,\mathrm{id}_{\palg}\oplus\ad\mathfrak
z(\kalg)|_{\palg}$. In Section~\ref{sec:parameter-families}, we prove
the existence for all parameters using Kac's realization of twisted
loop algebras, and construct the parameter isomorphisms from scalar
evaluations of the coproduct and spectral translations. The identity
involution is treated along the same lines in
Section~\ref{sec:identity-involution}. In Section~\ref{sec:minimal},
we prove the finite presentation, using a separation of positive
Cartan weights. Finally, the quantum duality conditions are verified
in Appendix~\ref{sec:quantum-duality-assumption}.

\medskip

\noindent {\bf Acknowledgements.}
The author used AI 
in preparing this paper, for exploring and checking arguments, editing the
exposition, and LaTeX typesetting. The author takes full responsibility for
the content of the paper and the correctness of its results. The author thanks Hongjia Chen, Ming Lu, and Weiqiang Wang for helpful comments and
suggestions. 

\section{Preliminaries}\label{sec:conventions}

\subsection{Enveloping algebras and relative primitives}
\label{ssec:basic-algebra}

Throughout the paper, we work over the field $\C$ of complex numbers.
For a
(possibly infinite-dimensional) complex Lie algebra $\mathfrak a$, we
write $\rU(\mathfrak a)$ and
$\rS(\mathfrak a)$ for its universal enveloping algebra and its
symmetric algebra. The PBW filtration of $\rU(\mathfrak a)$ counts
the number of Lie algebra factors, and its associated graded algebra
is $\rS(\mathfrak a)$. Both are Hopf algebras in which the elements
of $\mathfrak a$ are primitive; we write $\Delta_0$ and
$\varepsilon_0$ for their coproducts and counits, and
\[
 \nabla(X)=\Delta_0(X)-X\otimes1-1\otimes X
\]
for the reduced coproduct. An element $p$ of a Hopf algebra is
\emph{primitive} if $\Delta(p)=p\otimes1+1\otimes p$. All algebra
homomorphisms in this paper preserve the unit.

The symmetrization map
\[
 \operatorname{sym}:\rS(\mathfrak a)\to\rU(\mathfrak a),\qquad
 \operatorname{sym}(v_1\cdots v_n)
 =\tfrac1{n!}\sum_{\sigma\in \mathfrak S_n}v_{\sigma(1)}\cdots v_{\sigma(n)},
\]
is a linear isomorphism by the PBW theorem, and it is an isomorphism
of coalgebras:
\[
 \Delta_0\circ\operatorname{sym}
 =(\operatorname{sym}\otimes\operatorname{sym})\circ\Delta_0 .
\]
It maps $\rS(\mathfrak l)$ onto $\rU(\mathfrak l)$ for every Lie
subalgebra $\mathfrak l\subseteq\mathfrak a$. Every element of
$\rS(\mathfrak a)$ is a polynomial in finitely many basis vectors, so
the arguments below apply equally to infinite-dimensional Lie
algebras.

\begin{lemma}\label{lem:relative-primitives}
Let $\mathfrak l\subseteq\mathfrak a$ be a complex Lie subalgebra,
not necessarily finite-dimensional.  If $X\in\rU(\mathfrak a)$ satisfies
\[
 \nabla(X)\in \rU(\mathfrak l)\otimes\rU(\mathfrak a),
\]
then $X\in\rU(\mathfrak l)+\mathfrak a$.
More precisely, for any vector-space decomposition
$\mathfrak a=\mathfrak l\oplus\mathfrak m$, there is a unique
decomposition $X=Y+z$ with $Y\in\rU(\mathfrak l)$ and $z\in\mathfrak m$.
\end{lemma}

\begin{proof}
Since $\rU(\mathfrak a)$ is cocommutative, $\nabla(X)$ is invariant
under the flip of the tensor factors. Hence
\[
 \nabla(X)\in
 \bigl(\rU(\mathfrak l)\otimes\rU(\mathfrak a)\bigr)
 \cap\bigl(\rU(\mathfrak a)\otimes\rU(\mathfrak l)\bigr)
 =\rU(\mathfrak l)\otimes\rU(\mathfrak l).
\]
As $\operatorname{sym}$ is a coalgebra isomorphism mapping
$\rS(\mathfrak l)$ onto $\rU(\mathfrak l)$, the element
$x=\operatorname{sym}^{-1}(X)\in\rS(\mathfrak a)$ satisfies
$\nabla(x)\in\rS(\mathfrak l)\otimes\rS(\mathfrak l)$. Choose bases
$\{\ell_i\}$ of $\mathfrak l$ and $\{m_j\}$ of $\mathfrak m$, and
write $x=\sum_{n\geqslant0}x_n$ with $x_n\in\rS^n(\mathfrak a)$.

For a monomial $v_1\cdots v_n$, the component of its coproduct in
$\rS^{n-1}(\mathfrak a)\otimes\rS^1(\mathfrak a)$ is
\[
 \sum_{j=1}^n v_1\cdots\widehat{v_j}\cdots v_n\otimes v_j .
\]
For $n\geqslant2$ this is a component of the reduced coproduct, and only
$x_n$ contributes to it. Hence the projection of $\nabla(x)$ to
$\rS^{n-1}(\mathfrak a)\otimes\mathfrak m$ is
\[
 \sum_j\frac{\partial x_n}{\partial m_j}\otimes m_j ,
\]
which must vanish. As the $m_j$ are linearly independent,
$\partial x_n/\partial m_j=0$ for all $j$ and all $n\geqslant2$. It
follows that $x_n$ does not involve the variables $m_j$, that is,
$x_n\in\rS^n(\mathfrak l)$ for $n\geqslant2$. The constant term lies in
$\rS(\mathfrak l)$, and the linear term splits uniquely according to
$\mathfrak a=\mathfrak l\oplus\mathfrak m$. Applying
$\operatorname{sym}$, we obtain $X=Y+z$. The uniqueness follows from
$\rU(\mathfrak l)\cap\mathfrak m=0$, which is again a consequence of
the PBW theorem.
\end{proof}

\begin{corollary}\label{cor:primitives}
The primitive elements of $\rU(\mathfrak a)$ are exactly the
elements of $\mathfrak a$.
\end{corollary}

\begin{proof}
By Lemma~\ref{lem:relative-primitives} with $\mathfrak l=0$, a
primitive element has the form $c1+z$ with $c\in\C$ and
$z\in\mathfrak a$. Since $\nabla(c1)=-c\,1\otimes1$, we have $c=0$.
\end{proof}

\subsection{Symmetric pairs}\label{ssec:symmetric-pairs}

Let $\g$ be a finite-dimensional complex simple Lie algebra and let
$\omega$ be an involution of $\g$, possibly the identity. Write its
eigenspace decomposition as
\[
  \g=\kalg\oplus\palg,
  \qquad \omega|_{\kalg}=1,
  \qquad \omega|_{\palg}=-1,
\]
and also $\g^{+}=\kalg$, $\g^-=\palg$. Thus $\kalg=\g$ and
$\palg=0$ when $\omega=\mathrm{id}$. The identity involution is
allowed in Section~\ref{sec:criterion},
Section~\ref{sec:identity-involution} and
Appendix~\ref{sec:quantum-duality-assumption}; in
Sections~\ref{sec:rigidity}, \ref{sec:parameter-families}
and~\ref{sec:minimal} we assume $\omega\ne\mathrm{id}$. Extend
$\omega$ to $\g[u]$ by
\[
  \check\omega(xu^r)=(-1)^r\omega(x)u^r,
  \qquad\text{so that}\qquad
  \Ltw=\kalg[u^2]\oplus u\palg[u^2]
  =\bigoplus_{r\geqslant0}\g^{(-1)^r}u^r.
\]
We give $u$ degree one and call the induced nonnegative grading on
$\rU(\g[u])$ and $\rU(\Ltw)$ the \emph{current degree}. Unless
explicitly indicated as PBW degree, \emph{degree} means current
degree and its homogeneous extension to quantizations.

\begin{lemma}\label{lem:bracket-generation}
For every nontrivial involution $\omega$ of $\g$,
\[
 [\kalg,\palg]=\palg,
 \qquad
 [\palg,\palg]=\kalg.
\]
Consequently, $\Ltw$ is generated as a Lie algebra by $\kalg$ and
$u\palg$.
\end{lemma}

\begin{proof}
We have $[\kalg,\kalg]\subseteq\kalg$,
$[\kalg,\palg]\subseteq\palg$ and $[\palg,\palg]\subseteq\kalg$. Set
$\mathfrak a=[\palg,\palg]\oplus[\kalg,\palg]$. The Jacobi identity
gives $[\kalg,[\palg,\palg]]\subseteq[\palg,\palg]$ and
$[\kalg,[\kalg,\palg]]\subseteq[\kalg,\palg]$, while
\[
 [\palg,[\palg,\palg]]\subseteq[\palg,\kalg]
 =[\kalg,\palg],\qquad
 [\palg,[\kalg,\palg]]\subseteq[\palg,\palg].
\]
Thus $\mathfrak a$ is an ideal of $\g$. If $\mathfrak a=0$, then
$[\g,\palg]=0$, contradicting $\palg\ne0$ and $\mathfrak z(\g)=0$.
Since $\g$ is simple, we have $\mathfrak a=\g$. Taking the
$\pm1$-eigenspaces of $\omega$, we obtain $[\palg,\palg]=\kalg$ and
$[\kalg,\palg]=\palg$.

Let $\mathfrak l$ be the Lie subalgebra generated by $\kalg$ and
$u\palg$. Then $u^2\kalg=[u\palg,u\palg]\subset\mathfrak l$, and if
$u^{2r}\kalg\subset\mathfrak l$ for some $r\geqslant1$, then
$u^{2r+1}\palg=[u^{2r}\kalg,u\palg]\subset\mathfrak l$ and
$u^{2r+2}\kalg=[u^{2r+1}\palg,u\palg]\subset\mathfrak l$. By
induction, $\mathfrak l$ contains every homogeneous component of
$\Ltw$.
\end{proof}

For $\omega=\mathrm{id}$, the same induction shows that
$\Ltw=\g[u^2]$ is generated by $\g$ and $u^2\g$, as
$[\g,\g]=\g$.

Fix Chevalley generators $\{e_i,f_i,h_i\}_{i\in \I}$ of $\g$, let
$(c_{ij})_{i,j\in\I}$ be the Cartan matrix, and let $d_i$ be the
relatively prime positive integers such that $(d_ic_{ij})$ is
symmetric. Normalize the invariant form on $\g$ by
$(\alpha_i,\alpha_i)=2d_i$, and put
\begin{equation}\label{eq:root-vector-normalization}
 x_i^+=\sqrt{d_i}\,e_i,
 \qquad x_i^-=\sqrt{d_i}\,f_i,
 \qquad \xi_i=d_i h_i=[x_i^+,x_i^-].
\end{equation}
Then $(x_i^+,x_i^-)=1$ and $[\xi_i,x_j^\pm]=\pm d_i c_{ij}x_j^\pm$.

Throughout, $\{x_r\}$ and $\{x^r\}$ denote dual bases of $\g$, and
\[
 \Omega_{\g}=\sum_r x_r\otimes x^r .
\]
The standard rational Lie bialgebra structure on $\g[u]$ is given by
$\delta(x)=0$ and
\begin{equation}\label{eq:ambient-cobracket}
 \delta(xu^r)=
 \sum_{a+b=r-1}[x\otimes1,\Omega_{\g}]u^a v^b
 \in(\g\otimes\g)[u,v]\cong\g[u]\otimes\g[u]
 \qquad(x\in\g,\ r\geqslant 1).
\end{equation}

The involution preserves the invariant form, so $\kalg$ and $\palg$
are orthogonal and the form is nondegenerate on each of them. Let
$\{v_r\}$ and $\{v^r\}$ be dual bases of $\kalg$, and set
\begin{equation}\label{eq:k-Casimir-convention}
 \Omega_{\kalg}=\sum_r v_r\otimes v^r,
 \qquad
 C_{\kalg}=\sum_r v_rv^r.
\end{equation}
Then $[x\otimes1+1\otimes x,\Omega_{\kalg}]=0$ and
$[x,C_{\kalg}]=0$ for $x\in\kalg$, whereas for $y\in\palg$ the
commutator $[y,C_{\kalg}]$ is in general a nonzero quadratic element
of $\rU(\g)$.

For $y\in\palg$ we put
\begin{equation}\label{eq:mu-1-co}
 \gamma(yu):=-[1\otimes y,\Omega_{\kalg}]
 =-\sum_r v_r\otimes[y,v^r]\in\kalg\otimes\palg .
\end{equation}
By Lemma~\ref{lem:lift} below, this tensor is the first-order term of
the right coaction on a degree-one lift of $yu$. It is the value on
the generators $u\palg$ of the right coideal version of the Lie
coideal structure computed in \cite[\S4]{BR17}; we only need the
values \eqref{eq:mu-1-co}. The left coideal convention of
\cite[\S4]{BR17} uses $\tilde\gamma(yu)=[y\otimes1,\Omega_{\kalg}]$,
so that $\gamma(yu)=-\operatorname{flip}(\tilde\gamma(yu))$.

\subsection{The homogeneous Yangian}\label{ssec:yangian}

In this subsection, we recall the homogeneous Yangian, following the
reduced $\mathcal J$-presentation of \cite[\S~5]{GWX25}. The Yangian
and its classical limit are fixed as part of the quantization
problem. Let $\mathfrak h\subset\g$ and
$\Phi^+$ be the Cartan subalgebra and the set of positive roots.
Choose root vectors $x_\alpha^\pm\in\g_{\pm\alpha}$ with
$(x_\alpha^+,x_\alpha^-)=1$, and set
\begin{equation}\label{eq:nu}
 \nu(t)=\tfrac12\sum_{\alpha\in\Phi^+}
              \alpha(t)x_\alpha^-x_\alpha^+\in \rU(\g)
 \qquad(t\in\mathfrak h).    
\end{equation}
The \emph{homogeneous Yangian} $\Y$ is the unital
$\C[\hbar]$-algebra generated by two linear copies
$x\mapsto x$ and $x\mapsto \mathcal J(x)$ of $\g$, subject to
\[
 xy-yx=[x,y],\qquad
 [x,\mathcal J(y)]=\mathcal J([x,y])\qquad(x,y\in\g),
\]
and the following additional relations:
\begin{itemize}[leftmargin=2em]
\item $[\mathcal J(\xi_i),\mathcal J(\xi_j)]
 =\hbar^2[\nu(\xi_j),\nu(\xi_i)]$ for $i,j\in \I$, if $\g\not\cong\mathfrak{sl}_2$;
\item  $[[\mathcal J(e),\mathcal J(f)],\mathcal J(h)]
 =\hbar^2\bigl(f\mathcal J(e)-\mathcal J(f)e\bigr)h$, if
 $\g\cong\mathfrak{sl}_2$, where $e,f,h$ are the standard elements
 of $\mathfrak{sl}_2$.
\end{itemize}
It is graded by $\deg x=0$ and $\deg \mathcal J(x)=\deg\hbar=1$, and it
is a Hopf algebra with
\begin{align}
 &\Delta(x)=x\otimes1+1\otimes x, \qquad
 \varepsilon(x)=\varepsilon(\mathcal J(x))=0,\label{eq:couint}\\
 &\Delta(\mathcal J(x))
 =\mathcal J(x)\otimes1+1\otimes \mathcal J(x)
   +\tfrac12\hbar[x\otimes1,\Omega_{\g}].\label{eq:coproduct-J}
\end{align}
The existence and PBW results reviewed in \cite[\S~5.3]{GWX25}
(see also \cite{Lev93PBW}) show that $\Y$ is free over $\C[\hbar]$ and has classical limit
\[
 \cl:\Y/\hbar\Y\xrightarrow{\sim}\rU(\g[u]),
 \qquad \overline x\longmapsto x,
 \qquad \overline{\mathcal J(x)}\longmapsto xu.
\]
Thus $\mathcal J(x)$ is a fixed degree-one lift of $xu$. Here $u$ is
the current variable and $\hbar$ is a central deformation parameter.
Specializing $\hbar=1$ gives the ordinary
Yangian $\mathrm Y(\g)=\Y/(\hbar-1)\Y$. For $X\in\Y$ we write
\[
 \overline X=\cl(X\bmod\hbar)\in\rU(\g[u]),\qquad
 \nabla(X)=\Delta(X)-X\otimes1-1\otimes X ,
\]
and $\Delta^{\mathrm{op}}=\operatorname{flip}\circ\Delta$.

Every automorphism $\phi$ of $\g$ extends uniquely to a graded Hopf
algebra automorphism of $\Y$, again denoted $\phi$, with
$\phi(\mathcal J(x))=\mathcal J(\phi(x))$ and $\phi(\hbar)=\hbar$.
Indeed, the reduced presentation above is equivalent to Drinfeld's
$\mathcal J$-presentation \cite{Dri85}; see \cite[\S5]{GWX25}. The
defining relations of the latter only involve the Lie bracket and dual
bases for the invariant form, and hence they are preserved by $\phi$.
Moreover, \eqref{eq:coproduct-J} is compatible with $\phi$, since
$\Omega_\g$ is $(\phi\otimes\phi)$-invariant. In particular, $\omega$
acts on $\Y$. Note that on the classical level this action is
coefficient-wise and fixes $u$, so it is different from
$\check\omega$, which also changes the sign of $u$.

For $c\in\C$, the \emph{spectral translation}
\begin{equation}\label{eq:shift}
 \sigma_c(x)=x,\qquad \sigma_c(\mathcal J(x))=\mathcal J(x)+c\hbar x
 \qquad(x\in\g)
\end{equation}
is a graded Hopf algebra automorphism of $\Y$ with inverse
$\sigma_{-c}$, and its reduction modulo $\hbar$ is the identity of
$\rU(\g[u])$.

\subsection{Graded modules, saturation and lifting}
\label{ssec:graded}

Let $M=\bigoplus_{d\geqslant 0}M_d$ be a graded module over
$\C[\hbar]$, so that multiplication by $\hbar$ maps $M_{d-1}$ to
$M_d$. We write $\overline m$ for the class of $m$ modulo $\hbar M$
and $\overline f$ for an induced map. For homogeneous $m\in M_d$,
\[
 \overline m=0\quad\Longleftrightarrow\quad
 m=\hbar m_1\quad\text{ for some }m_1\in M_{d-1},
\]
and $m_1$ is unique if multiplication by $\hbar$ is injective. Since
each division lowers the degree, repeated division stops after at
most $d+1$ steps. We call this argument \emph{degree descent}; note
that it does not require the homogeneous components to be
finite-dimensional. A $\C[\hbar]$-module is
torsion-free if multiplication by every nonzero polynomial in $\hbar$
is injective; free modules are torsion-free.

Let $M$ be a $\C[\hbar]$-module on which multiplication by $\hbar$ is
injective, and let $N\subseteq M$ be a submodule. We say that $N$ is
\emph{$\hbar$-saturated} (or simply \emph{saturated}) in $M$ if
\[
 N\cap\hbar M=\hbar N,
\]
which is equivalent to injectivity of the natural map
$N/\hbar N\to M/\hbar M$.

\begin{lemma}\label{lem:tensor-saturation}
Let $N\subseteq M$ be saturated, where $\hbar$ acts injectively on
$M$, and let $P$ be a free $\C[\hbar]$-module. Then
$N\otimes_{\C[\hbar]}P\to M\otimes_{\C[\hbar]}P$ is injective, and
its image is saturated in $M\otimes_{\C[\hbar]}P$.
\end{lemma}

\begin{proof}
A $\C[\hbar]$-basis of $P$ identifies the map with a direct sum of
copies of the inclusion $N\subseteq M$.
\end{proof}

The following lemma lifts an isomorphism at the classical level to an
isomorphism before reduction. The argument terminates because the
grading is nonnegative, and the torsion-freeness of the target is used
to cancel $\hbar$ in the proof of injectivity.

\begin{lemma}\label{lem:graded-lifting}
Let $f:A\to B$ be a graded homomorphism of nonnegatively graded
$\C[\hbar]$-algebras, where $\deg\hbar=1$, and suppose that $B$ is
$\C[\hbar]$-torsion-free.  If $\bar f:A/\hbar A\to B/\hbar B$ is an
isomorphism, then $f$ is an isomorphism.  The same statement holds
for a graded $\C[\hbar]$-linear map of nonnegatively graded
$\C[\hbar]$-modules.
\end{lemma}

\begin{proof}
A bijective algebra homomorphism is an algebra isomorphism, so it
suffices to prove the module statement. Write $A_n,B_n$ for the
homogeneous components. In degree zero,
$(\hbar A)_0=(\hbar B)_0=0$, so $f:A_0\to B_0$ is an isomorphism.

For surjectivity, suppose that every element of $B_{n-1}$ lifts to
$A_{n-1}$ and take $b\in B_n$. Choose $a_0\in A_n$ with
$\overline{f(a_0)}=\bar b$. Then $b-f(a_0)=\hbar b_1$ with
$b_1\in B_{n-1}$, and by induction $b_1=f(a_1)$ for some
$a_1\in A_{n-1}$. Hence $b=f(a_0+\hbar a_1)$.

For injectivity, let $a\in A_n$ satisfy $f(a)=0$. Since $\bar f$ is
injective, we have $a=\hbar a_1$ with $a_1\in A_{n-1}$. As
multiplication by $\hbar$ is injective on $B$, the equality
$\hbar f(a_1)=0$ implies $f(a_1)=0$. By induction, $a_1=0$, and thus
$a=0$. Because every element is a finite sum of homogeneous elements,
this proves the lemma.
\end{proof}

\subsection{Coideal quantizations}\label{ssec:coideal-quantizations}

Write $\inc:\rU(\Ltw)\hookrightarrow \rU(\g[u])$ for the natural
inclusion. We also write $\Delta_0$ for its restriction
$\rU(\Ltw)\to \rU(\Ltw)\otimes \rU(\g[u])$ and call it the
\emph{classical coaction}.

\begin{definition}
\label{def:polynomial-coq}
A \emph{homogeneous coideal quantization} of $\Ltw$ in $\Y$, or a
\emph{twisted Yangian} associated with $(\g,\kalg)$, is a graded
unital subalgebra $B\subset\Y$ such that
\begin{enumerate}[label=\textup{(\roman*)}]
\item $B$ is a right coideal: $\Delta(B)\subset B\otimes_{\C[\hbar]}\Y$;
\item $B$ is free as a graded $\C[\hbar]$-module;
\item the natural reduction map identifies
\begin{equation}\label{eq:coideal-classical-inclusion}
 B/\hbar B\ \cong\ \rU(\Ltw) \subset
 \rU(\g[u])\cong\Y/\hbar\Y,
\end{equation}
that is, $B/\hbar B\to\Y/\hbar\Y$ is injective and $\cl$ maps its
image onto $\rU(\Ltw)$.
\end{enumerate}
\end{definition}

Note that the injectivity in (iii) is equivalent to the saturation
condition
\begin{equation}\label{eq:coideal-saturation}
 B\cap\hbar\Y=\hbar B.
\end{equation}

\begin{remark}\label{rem:definition}
This is the right coideal version of \cite[Def.~3.15]{LWW25}
for the fixed Yangian $\Y$. Note that condition
\eqref{eq:coideal-saturation} depends on the embedding, not only on
the abstract module. For instance, $\hbar\C[\hbar]\subset\C[\hbar]$ is
a free submodule, but the nonzero class of $\hbar$ in
$\hbar\C[\hbar]/\hbar^2\C[\hbar]$ maps to zero in
$\C[\hbar]/\hbar\C[\hbar]$. Conversely, in our nonnegatively graded
embedded setting, (ii) follows from (iii). Indeed, lifting an ordered
homogeneous PBW basis of $\rU(\Ltw)$, we obtain a graded map to $B$
from a free module which is an isomorphism modulo $\hbar$, and hence
an isomorphism by Lemma~\ref{lem:graded-lifting}. We keep (ii) to be
consistent with \cite{LWW25}. As in \cite{LWW25}, neither a preferred
degree-one lift nor a formula for the first-order coaction is part of
the definition. In contrast, \cite[Def.~5.3]{BR17} additionally
imposes a condition on the coproduct modulo $\hbar^2$, which
prescribes the first-order coaction by the Lie coideal structure; see
\cite[Rem.~3.17]{LWW25}. By Lemma~\ref{lem:lift} and
Remark~\ref{rem:semiclassical}, in the right coideal convention, this
condition holds automatically for suitable lifts.
\end{remark}

For a homogeneous coideal quantization $B$, we write
$\cli:B/\hbar B\xrightarrow{\sim}\rU(\Ltw)$ for the identification
\eqref{eq:coideal-classical-inclusion} and $\rho=\Delta|_B:B\to
B\otimes_{\C[\hbar]}\Y$ for its coaction. The coassociativity and
counit identities
\[
 (\rho\otimes\mathrm{id})\rho=(\mathrm{id}\otimes\Delta)\rho,
 \qquad (\mathrm{id}\otimes\varepsilon)\rho=\mathrm{id}_B
\]
follow from those of $\Y$; see \cite[Rem.~3.16]{LWW25}. The
reduction of $\rho$ is the classical coaction $\Delta_0$. Since
$(\hbar B)_0=0$, reduction is injective in degree zero, and hence
\[
 B_0=\rU(\kalg),\qquad \Y_0=\rU(\g).
\]
As $\Y$ is free, the reduction of $B\otimes_{\C[\hbar]}\Y$ is
$(B/\hbar B)\otimes_{\C}(\Y/\hbar\Y)\cong\rU(\Ltw)\otimes\rU(\g[u])$.
Here and below, tensor products are taken over $\C[\hbar]$ before
reduction and over $\C$ after reduction.

\begin{remark}[Abstract realizations]\label{rem:realizations}
When comparing presentations, it is convenient to consider an
abstract graded algebra $B$ together with an injective graded homomorphism
$\varphi:B\to\Y$ whose image satisfies
Definition~\ref{def:polynomial-coq}. The maps $\cli$ and $\rho$ are
then characterized by $\cl\circ\overline\varphi=\inc\circ\cli$ and
$(\varphi\otimes\mathrm{id})\circ\rho=\Delta\circ\varphi$. An
\emph{embedding-preserving isomorphism} between two realizations
$(B_1,\varphi_1)$ and $(B_2,\varphi_2)$ is a graded algebra isomorphism
$F:B_1\to B_2$ with $\varphi_2\circ F=\varphi_1$. Then
$\varphi_1(B_1)=\varphi_2(B_2)$, and $F$ automatically satisfies
\[
 \mathrm{cl}_{\imath,2}\circ\overline F=\mathrm{cl}_{\imath,1},\qquad \rho_2\circ F=(F\otimes\mathrm{id})\circ\rho_1.
\]
Indeed, reducing $\varphi_2F=\varphi_1$ modulo $\hbar$ gives
$\inc\circ\mathrm{cl}_{\imath,2}\circ\overline F=\inc\circ\mathrm{cl}_{\imath,1}$, and $\inc$ is
injective. Moreover
$(\varphi_2\otimes\mathrm{id})\rho_2F=\Delta\varphi_2F=\Delta\varphi_1
=(\varphi_2\otimes\mathrm{id})(F\otimes\mathrm{id})\rho_1$, and
$\varphi_2\otimes\mathrm{id}_{\Y}$ is injective because $\Y$ is free
over $\C[\hbar]$. An isomorphism of right $\Y$-comodule algebras is
a graded algebra isomorphism $F$ with
$\rho_2\circ F=(F\otimes\mathrm{id})\circ\rho_1$, where the Yangian
is fixed. In contrast, \emph{uniqueness up to graded algebra
isomorphism} only compares the underlying graded $\C[\hbar]$-algebras,
and the embeddings and the induced classical identifications are
allowed to vary. Outside this remark, we identify a twisted Yangian
with its image in $\Y$.
\end{remark}

\begin{remark}[Filtered twisted Yangians]\label{rem:filtered}
Let $\mathrm Y(\g)=\Y/(\hbar-1)\Y$ be the ordinary Yangian with its
standard filtration, of which $\Y$ is the Rees algebra. Homogeneous
coideal quantizations correspond, via specialization at $\hbar=1$ and
the Rees construction, to filtered right coideal subalgebras of
$\mathrm Y(\g)$ with associated graded algebra
$\rU(\Ltw)\subset\rU(\g[u])$, where the saturation identifies the two
filtrations. Thus all results of this paper apply to filtered twisted
Yangians.
\end{remark}

\section{Affine Serre relations and a saturation criterion}
\label{sec:criterion}

In this section, $\omega$ is an arbitrary involution of $\g$, possibly
the identity. We prove a criterion for $\hbar$-saturation and the PBW
property, which will be applied to all twisted Yangians in this
paper: to nontrivial involutions in
Section~\ref{sec:parameter-families} and to the identity involution in
Section~\ref{sec:identity-involution}. The classical part consists of
generators of the relation ideal of $\Ltw$ and their weights; see
\S\ref{ssec:free-model}--\S\ref{ssec:classical-relations}. In the
quantum part, we use these weights to exclude a first failure of
saturation; see \S\ref{ssec:saturation-pbw}. Note that we do not need
any quantum relations among the generators. We only use the classical
affine Serre relations, the coproduct and degree descent, and
Proposition~\ref{prop:saturation-criterion}(d) shows a posteriori that
suitable lifts of the Serre relations are defining relations. For this
reason, the criterion can be applied to an arbitrary coideal
subalgebra in Sections~\ref{sec:rigidity}
and~\ref{sec:identity-involution}; compare \cite[Rem.~5.8]{BR17},
where flatness is deduced from explicit quantum relations.

\subsection{The free model}\label{ssec:free-model}

Put
\[
 (e,V)=(1,\palg)\quad\text{if }\omega\ne\mathrm{id},\qquad
 (e,V)=(2,\g)\quad\text{if }\omega=\mathrm{id}.
\]
In both cases $\Ltw_0=\kalg$, $\Ltw_d=0$ for $0<d<e$, and
$\Ltw_e=u^eV$. Moreover, by Lemma~\ref{lem:bracket-generation} and the
remark following it, $\Ltw$ is generated as a Lie algebra by $\kalg$
and $u^eV$.

Let $\underline V=\{\underline v:v\in V\}$ be a copy of $V$ placed in
degree $e$, and let
\[
 \mathfrak F=\kalg\ltimes\operatorname{Lie}(\underline V)
\]
be the semidirect product of $\kalg$ with the free Lie algebra on
$\underline V$, where $\kalg$ acts on $\underline V$ through its
adjoint action on $V$. Thus $\mathfrak F$ is generated by $\kalg$ and
$\underline V$, subject to the relations of $\kalg$ and to
\[
 [x,\underline v]=\underline{[x,v]}\qquad(x\in\kalg,\ v\in V).
\]
A Lie word has degree $e$ times the number of its letters from
$\underline V$, so that $\mathfrak F_0=\kalg$, $\mathfrak F_e=\underline V$
and $\mathfrak F_d=0$ for $0<d<e$. The evaluation map
\[
 \pi_0:\mathfrak F\longrightarrow\Ltw,
 \qquad \pi_0(x)=x,\quad \pi_0(\underline v)=u^ev,
\]
is a well-defined graded Lie algebra homomorphism, because the defining
relations of $\mathfrak F$ hold in $\Ltw$. It is surjective, and it
is bijective in degrees at most $e$. Its kernel consists of the
additional relations of $\Ltw$.

\subsection{Affine data and Serre words}
\label{ssec:classical-relations}

Let
\[
 L=\g[u,u^{-1}]^{\check\omega}
 =\bigoplus_{m\in\mathbb Z}\g^{(-1)^m}u^m
\]
be the twisted loop algebra, graded by the power of $u$. Its
nonnegative part is $\Ltw$. We describe the relation ideal of $\Ltw$
using the following data.

\begin{definition}\label{def:affine-datum}
An \emph{affine datum} for $(\g,\omega)$ consists of an affine
Kac--Moody algebra $\widehat\g$ with Chevalley generators
$e_i,f_i,h_i$ ($i\in\hI$), generalized Cartan matrix
$A=(a_{ij})_{i,j\in\hI}$ with $[h_i,e_j]=a_{ij}e_j$, derived algebra
$\widehat\g'$ generated by the $e_i,f_i$, and central element
$\mathsf c$ spanning the center of $\widehat\g'$, together with a
Lie algebra isomorphism
\[
 \lambda:\widehat\g'/\C\mathsf c\xrightarrow{\ \sim\ }L,
\]
such that, writing $\overline e_i,\overline f_i,\overline h_i$ for
the images of $e_i,f_i,h_i$, the following hold.
\begin{enumerate}[label=\textup{(A\arabic*)}]
\item Each $\overline e_i$ is homogeneous of degree $s_i\in\{0,e\}$,
each $\overline f_i$ of degree $-s_i$, and each $\overline h_i$ of
degree $0$.
\item The Lie algebra $\kalg$ is generated by
$\mathfrak t=\operatorname{span}_\C\{\overline h_i:i\in\hI\}$ and by
$\overline e_j,\overline f_j$ for $j\in\hI^0=\{j:s_j=0\}$.
\item For $i\in\hI^e=\{i:s_i=e\}$ write $\overline e_i=u^ey_i$ with
$y_i\in V$. Then $V$ is spanned by the iterated brackets
\[
 [x_1,[x_2,\ldots,[x_k,y_i]\cdots]]\qquad
 \bigl(k\geqslant0,\ i\in\hI^e,\ x_1,\ldots,x_k\in\{\overline e_j:j\in\hI^0\}\bigr).
\]
\end{enumerate}
\end{definition}

Every involution $\omega$ admits an affine datum. For
$\omega\ne\mathrm{id}$, this follows from Kac's classification of
finite-order automorphisms and the associated loop realizations; see
\cite[Thms.~8.6 and~8.7]{Kac90} and \cite[Prop.~8.6]{Kac90}.
Conditions (A1)--(A3) are verified in
Lemma~\ref{lem:affine-degree-one}. For $\omega=\mathrm{id}$, the
standard untwisted loop realization, with loop variable $v=u^2$,
gives an affine datum with $e=2$; see
Lemma~\ref{lem:identity-classical-relations}. In either case, we fix
such an affine datum for the rest of this section.

By (A1), we have $\mathfrak t\subset L_0=\kalg$. The
abelian subalgebra $\mathfrak t$ acts semisimply on $L$, because
$\widehat\g'$ is the sum of its weight spaces for
$\mathfrak h'=\operatorname{span}_\C\{h_i\}$; hence it acts
semisimply on $\g$ and on $\Ltw$. For $\mu\in\mathfrak t^*$ and
$\epsilon\in\{+1,-1\}$ we write
\[
 (\g^\epsilon)_\mu
 =\{x\in\g^\epsilon:[t,x]=\mu(t)x\text{ for all }t\in\mathfrak t\}.
\]
In this section, weights always mean $\mathfrak t$-weights.

To distinguish formal words from their values, we set
\[
 E_i=
 \begin{cases}
 \overline e_i\in\kalg=\mathfrak F_0,&i\in\hI^0,\\
 \underline{y_i}\in\mathfrak F_e,&i\in \hI^e,
 \end{cases}
\]
so that $\pi_0(E_i)=\overline e_i$ in both cases. For $i\ne j$, define
the \emph{Serre word}
\begin{equation}\label{eq:classical-serre-words-pbw}
 r_{ij}=(\ad E_i)^{1-a_{ij}}E_j\in\mathfrak F,\qquad
 d_{ij}:=\deg r_{ij}=(1-a_{ij})s_i+s_j.
\end{equation}
Consider the finite index set
\begin{equation}\label{eq:serre-index-set}
 \mathcal S=\{(i,j)\in\hI\times\hI:\ i\ne j,\ d_{ij}>e\},
\end{equation}
and for $\nu=(i,j)\in\mathcal S$ write $r_\nu=r_{ij}$ and
$d_\nu=d_{ij}$.

\begin{lemma}
\label{lem:classical-relations}
The Serre words $r_\nu$, $\nu\in\mathcal S$, generate $\ker\pi_0$ as a
Lie ideal of $\mathfrak F$, and generate the kernel of
$\rU(\mathfrak F)\to\rU(\Ltw)$ as a two-sided ideal of
$\rU(\mathfrak F)$.
\end{lemma}

\begin{proof}
The Serre relations of $\widehat\g'$, and hence of $L$, give
\[
 \pi_0(r_{ij})=(\ad\overline e_i)^{1-a_{ij}}\overline e_j=0,
\]
so all Serre words lie in $\ker\pi_0$. Since $s_i\in\{0,e\}$, each
$d_{ij}$ is a multiple of $e$. If $d_{ij}\leqslant e$, then $r_{ij}$
lies in $\mathfrak F_0\oplus\mathfrak F_e$, where $\pi_0$ is
injective, and thus $r_{ij}=0$ already in $\mathfrak F$. In other
words, these relations already follow from the relations of $\kalg$
and its action on $V$. Therefore, we only need to consider the words
$r_\nu$ with $\nu\in\mathcal S$.

Let $\mathfrak I$ be the Lie ideal of $\mathfrak F$ generated by the
$r_\nu$, $\nu\in\mathcal S$. It is spanned by the $r_\nu$ and their
iterated brackets with elements of $\mathfrak F$. Put
$\mathfrak A=\mathfrak F/\mathfrak I$, and let
$\overline\pi_0:\mathfrak A\to\Ltw$ be the induced graded
surjection. We show that $\overline\pi_0$ is injective.

\emph{Step 1: the positive part.}
Let $\widehat{\mathfrak n}_+$ be the subalgebra of $\widehat{\g}'$
generated by the $e_i$. By the Gabber--Kac theorem
\cite[Thm.~9.11]{Kac90}, its defining relations are the Serre
relations $(\ad e_i)^{1-a_{ij}}e_j=0$ for $i\ne j$. These relations
hold for the elements $E_i+\mathfrak I$ of $\mathfrak A$. Indeed,
those of degree at most $e$ already hold in $\mathfrak F$, and the
remaining ones are the generators of $\mathfrak I$. Thus
$e_i\mapsto E_i+\mathfrak I$ induces a homomorphism
\[
 \alpha:\widehat{\mathfrak n}_+\longrightarrow\mathfrak A,
 \qquad \mathcal N:=\alpha(\widehat{\mathfrak n}_+).
\]
The composite $\overline\pi_0\circ\alpha$ and the restriction of
$\widehat\g'\to\widehat\g'/\C\mathsf c\xrightarrow{\lambda}L$ to
$\widehat{\mathfrak n}_+$ agree on the generators $e_i$, and hence on
$\widehat{\mathfrak n}_+$. Here the latter takes values in
$\Ltw$, as $s_i\geqslant0$. By the triangular decomposition of
$\widehat\g'$, we have $\widehat{\mathfrak n}_+\cap\C\mathsf c=0$, so
this restriction is injective. Therefore, $\overline\pi_0$ is
injective on $\mathcal N$.

\emph{Step 2: $\mathfrak A=\kalg+\mathcal N$.}
Here $\kalg$ denotes its image in $\mathfrak A$. By the relation
$[x,\underline v]=\underline{[x,v]}$ and (A3), every element of
$\underline V$ is a linear combination of iterated brackets of the
$E_j$, $j\in\hI^0$, applied to the $E_i$, $i\in\hI^e$. Hence
$\underline V\subset\mathcal N$. It remains to show that
$\kalg+\mathcal N$ is closed under brackets. Clearly
$[\kalg,\kalg]\subseteq\kalg$ and
$[\mathcal N,\mathcal N]\subseteq\mathcal N$. By the Jacobi identity,
the set of $x\in\kalg$ with $[x,\mathcal N]\subseteq\kalg+\mathcal N$
is a Lie subalgebra of $\kalg$. Thus by (A2), it suffices to check
the generators. The elements of $\mathfrak t$ act on each $E_i$ by
scalars, so they preserve $\mathcal N$. We also have
$\overline e_j=E_j\in\mathcal N$ for $j\in\hI^0$. For $j\in\hI^0$,
the relations $[f_j,e_i]=-\delta_{ij}h_j$ in $\widehat\g'$ give
\[
 [\overline f_j,E_i]=-\delta_{ji}\overline h_j\quad(i\in\hI^0),
 \qquad
 [\overline f_j,E_i]=\underline{[\overline f_j,y_i]}=0\quad(i\in\hI^e).
\]
By induction on the length of Lie words in the $E_i$ and the Jacobi
identity, we obtain $[\overline f_j,\mathcal N]\subseteq\mathcal
N+\C\overline h_j$. Thus $\kalg+\mathcal N$ is a Lie subalgebra of
$\mathfrak A$ containing both generating subspaces, and therefore
$\kalg+\mathcal N=\mathfrak A$.

\emph{Step 3: conclusion.}
In degree zero, $\mathfrak A_0=\kalg$ and $\overline\pi_0$ is the
identity. For $d>0$, Step~2 gives $\mathfrak A_d=\mathcal N_d$, on
which $\overline\pi_0$ is injective by Step~1. Hence
\[
 \mathfrak F/\mathfrak I\cong\Ltw,
 \qquad \ker\pi_0=\mathfrak I.
\]
For the enveloping algebras, note that a homomorphism from
$\rU(\mathfrak F)$ that kills the $r_\nu$ also kills their iterated
commutators, and thus it factors through
$\rU(\mathfrak F/\mathfrak I)=\rU(\Ltw)$. Therefore,
$\rU(\Ltw)\cong\rU(\mathfrak F)/\langle r_\nu\rangle$, where
$\langle r_\nu\rangle$ denotes the two-sided ideal generated by the
$r_\nu$.
\end{proof}

Next, we determine the weights of the Serre words $r_\nu$.

\begin{lemma}
\label{lem:classical-relation-weights}
Each $r_\nu$ is a $\mathfrak t$-weight vector, and its weight
$\mu_\nu$ and degree $d_\nu$ satisfy
\begin{equation}\label{eq:classical-serre-weight-exclusion}
\mu_\nu\ne 0,\qquad  \bigl(\g^{(-1)^{d_\nu}}\bigr)_{\mu_\nu}=0.
\end{equation}
\end{lemma}

\begin{proof}
Let $Q=\bigoplus_{i\in\hI}\mathbb Z\alpha_i$ be the root lattice of
$\widehat\g$. The algebra $\widehat\g'$ is $Q$-graded, with $e_i$,
$f_i$, $h_i$ in degrees $\alpha_i$, $-\alpha_i$, $0$; its component
of degree $\beta\ne0$ is the root space $\widehat\g_\beta$, which is
zero unless $\beta$ is a root. For $\beta=\sum_jk_j\alpha_j\in Q$ put
\[
 \beta(h_i)=\sum_jk_ja_{ij},\qquad
 \deg_s\beta=\sum_jk_js_j .
\]
Then $h_i$ acts on the $\beta$-component by $\beta(h_i)$, and by
(A1) the map $\lambda$ sends the $\beta$-component into the component
of $L$ of degree $\deg_s\beta$ and $\mathfrak t$-weight
$\bar\beta$, where $\bar\beta(\overline h_i)=\beta(h_i)$. The weight
$\bar\beta$ is well defined. Indeed, if $\sum_ic_i\overline h_i=0$,
then
$\sum_ic_ih_i\in\C\mathsf c$, and $\beta(\mathsf c)=0$ for all
$\beta\in Q$.

We claim that $\beta\in Q$ is determined by the pair
$(\bar\beta,\deg_s\beta)$. Indeed, $\bar\beta=0$ means
$\sum_jk_ja_{ij}=0$ for all $i$. For an affine Cartan matrix the
integral solutions of this system are the integer multiples of the
vector of marks $(m_j)$, whose entries are positive integers
\cite[Ch.~4]{Kac90}; equivalently, $\beta\in\mathbb Z\delta$ with
$\delta=\sum_jm_j\alpha_j$. Some $s_j$ equals $e>0$, since otherwise
$\lambda(\widehat\g')$ would lie in $L_0$. Hence
$\deg_s\delta=\sum_jm_js_j>0$, and $\beta=0$ if moreover
$\deg_s\beta=0$. Applying this to the difference of two elements of
$Q$, we obtain the claim.

Each $E_i$ is a weight vector of weight $\mu_i=\bar\alpha_i$ and
degree $s_i$. Therefore, $r_{ij}$ is a weight vector of weight
\[
 \mu_{ij}=(1-a_{ij})\mu_i+\mu_j=\bar\beta_{ij},\qquad
 \beta_{ij}=(1-a_{ij})\alpha_i+\alpha_j,
\]
and $d_{ij}=\deg_s\beta_{ij}$. Note that $\beta_{ij}$ is not a root.
Indeed, since $\alpha_j-\alpha_i$ is not a root (its coefficients have
opposite signs), the root-string property for the simple real root
$\alpha_i$ shows that the roots of the form $\alpha_j+k\alpha_i$ are
exactly $\alpha_j,\alpha_j+\alpha_i,\ldots,\alpha_j-a_{ij}\alpha_i$.
The weight $\beta_{ij}$ is one step beyond the end of this string.
Moreover, $\beta_{ij}\ne0$.

Suppose now that $0\ne x\in(\g^{(-1)^{d}})_{\mu}$, where
$(\mu,d)=(\mu_{ij},d_{ij})$. Then $xu^d$ is a nonzero element of $L$
of weight $\mu$ and degree $d$. Choose $\tilde x\in\widehat\g'$ with
$\lambda(\tilde x+\C\mathsf c)=xu^d$ and decompose $\tilde x$ into
$Q$-homogeneous components. As $L$ is the direct sum of its
components of fixed degree and weight, $xu^d$ is the image of the sum
of the components of degrees $\beta$ with $\bar\beta=\mu$ and
$\deg_s\beta=d$. By the claim, $\beta_{ij}$ is the only such element
of $Q$, and $\widehat\g'_{\beta_{ij}}=0$. This is a contradiction,
which proves the second part of
\eqref{eq:classical-serre-weight-exclusion}.

Finally, $\overline h_i\in\mathfrak t$ and
\begin{equation}\label{eq:serre-nonzero-cartan-weight}
 \mu_{ij}(\overline h_i)=\beta_{ij}(h_i)
   =2(1-a_{ij})+a_{ij}=2-a_{ij}\ne0,
\end{equation}
because $a_{ij}\leqslant0$ for $i\ne j$.
\end{proof}

Note that the weight exclusion
\eqref{eq:classical-serre-weight-exclusion} holds for the chosen
generators $r_\nu$ of the ideal, but not for every element of
$\ker\pi_0$; see Example~\ref{ex:sl2-serre}.

\subsection{The saturation criterion}\label{ssec:saturation-pbw}

Fix an affine datum for $(\g,\omega)$, with Serre words $r_\nu$ of
degrees $d_\nu$ and weights $\mu_\nu$. Let $\sfb:V\to\Y_e$ be a
$\kalg$-equivariant linear map, $[x,\sfb(v)]=\sfb([x,v])$, such that
for all $v\in V$:
\begin{enumerate}[label=\textup{(Q\arabic*)}]
\item $\overline{\sfb(v)}=u^ev$;
\item $\omega(\sfb(v))=(-1)^e\sfb(v)$, where $\omega$ acts on $\Y$ as
in \S\ref{ssec:yangian};
\item $\Delta(\sfb(v))-\sfb(v)\otimes1\in
\rU(\kalg)[\hbar]\otimes_{\C[\hbar]}\Y$.
\end{enumerate}
Let $\B\subset\Y$ be the $\C[\hbar]$-subalgebra generated by
$\rU(\kalg)$ and $\sfb(V)$. Clearly $\B$ is graded. By (Q3), $\B$ is a
right coideal subalgebra, since $\rU(\kalg)[\hbar]\subset\B$ and the
elements of $\kalg$ are primitive. Moreover, $\B$ is spanned by
monomials in $\mathfrak t$-weight vectors of $\kalg$ and of $\sfb(V)$,
and by (Q2), each of them is an $\omega$-eigenvector. Hence $\B$ is
the direct sum
of its intersections with the simultaneous components of $\Y$ of
fixed degree, $\mathfrak t$-weight and $\omega$-eigenvalue.

By the relations of $\kalg$ and the equivariance of $\sfb$, we have a
Lie algebra homomorphism $\mathfrak F\to\B$, where $\B$ is regarded as
a Lie algebra with the commutator bracket. By the universal property
of enveloping algebras, it extends to a graded surjection
\begin{equation}\label{eq:quantum-word-evaluation}
 \pi:\mathscr F:=\C[\hbar]\otimes\rU(\mathfrak F)
       \twoheadrightarrow \B,\qquad
 x\longmapsto x,\quad\underline v\longmapsto \sfb(v),\quad
 \hbar\longmapsto\hbar.
\end{equation}
By (Q1), the composite
\[
 \rU(\mathfrak F)
 \xrightarrow{\pi\bmod\hbar}\B/\hbar\B
 \longrightarrow\Y/\hbar\Y
 \xrightarrow{\cl}\rU(\g[u])
\]
coincides with the homomorphism induced by $\pi_0$, followed by
the natural inclusion $\rU(\Ltw)\hookrightarrow\rU(\g[u])$. In
particular, $\cl$ maps the image of $\B$ in $\Y/\hbar\Y$ onto
$\rU(\Ltw)$. Here the middle arrow is the natural map
$\B/\hbar\B\to\Y/\hbar\Y$ induced by the inclusion $\B\subset\Y$; we
will show that it is injective by the saturation argument below.
For each Serre word, we put $\mathscr R_\nu=\pi(r_\nu)$.
Then $\overline{\mathscr R_\nu}=0$, so $\mathscr R_\nu\in\hbar\Y$,
and $\mathscr R_\nu$ has degree $d_\nu$, weight $\mu_\nu$ and
$\omega$-eigenvalue $(-1)^{d_\nu}$.

To identify the image of $\B$ with $\B/\hbar\B$, we need to show that
$\B\cap\hbar\Y=\hbar\B$, that is, an element of $\B$ divisible by
$\hbar$ in $\Y$ is also divisible by $\hbar$ in $\B$. We first prove
two lemmas.

\begin{lemma}
\label{lem:low-degree-saturation}
Let $n\geqslant1$, and suppose that
$\B_j\cap\hbar\Y_{j-1}=\hbar \B_{j-1}$ for $0\leqslant j<n$ (for
$j=0$ this is automatic, as $\Y_{-1}=0$). Then the graded
$\C[\hbar]$-submodule
\[
 \mathcal M_n=\sum_{j=0}^{n-1}\C[\hbar]\B_j\subset\Y
\]
is saturated in $\Y$ and has homogeneous components
\begin{equation}\label{eq:low-degree-saturated-module}
 (\mathcal M_n)_j=
 \begin{cases}
 \B_j,&0\leqslant j<n,\\
 \hbar^{j-n+1}\B_{n-1},&j\geqslant n.
 \end{cases}
\end{equation}
Moreover,
\begin{equation}\label{eq:low-degree-module-coaction}
 \Delta(\mathcal M_n)\subset\mathcal M_n\otimes\Y,
\end{equation}
and $\mathcal M_n\otimes\Y$ is saturated in $\Y\otimes\Y$.
\end{lemma}

\begin{proof}
The formula \eqref{eq:low-degree-saturated-module} follows from
$\hbar \B_j\subset \B_{j+1}$. To prove the saturation, let
$\hbar v\in(\mathcal M_n)_j$ with $v\in\Y_{j-1}$. If $j<n$, then by
the assumption in degree $j$, we have
$v\in\B_{j-1}=(\mathcal M_n)_{j-1}$. If $j\geqslant n$, we write
$\hbar v=\hbar^{j-n+1}b$ with $b\in\B_{n-1}$ and cancel $\hbar$ in
$\Y$. Then $v=\hbar^{j-n}b\in(\mathcal M_n)_{j-1}$, which also covers
the case $j=n$. The case of nonhomogeneous elements follows by taking
homogeneous components.

For \eqref{eq:low-degree-module-coaction}, let $b\in\B_j$ with $j<n$.
By the graded coideal property, every homogeneous first factor of
$\Delta(b)$ lies in $\B$ and has degree at most $j$, and hence lies in
$\mathcal M_n$. Clearly this inclusion is preserved by multiplication
by $\hbar$. The last statement follows from
Lemma~\ref{lem:tensor-saturation}, as $\Y$ is free over
$\C[\hbar]$.
\end{proof}

\begin{lemma}
\label{lem:excluded-weight-descent}
Under the assumptions of Lemma~\ref{lem:low-degree-saturation},
let $\epsilon\in\{+1,-1\}$ and $\mu\in\mathfrak t^*$ satisfy
$(\g^\epsilon)_\mu=0$.
If $S\in\Y_m$, $0\leqslant m<n$, has weight $\mu$,
$\omega(S)=\epsilon S$, and
\[
 \nabla(S)\in\mathcal M_n\otimes\Y,
\]
then $S\in \B_m$.
\end{lemma}

\begin{proof}
We prove the lemma by induction on $m$. We first show that
$\overline S\in\rU(\Ltw)_m$. Reducing
$\nabla(S)\in\mathcal M_n\otimes\Y$ modulo $\hbar$, we obtain
\[
 \Delta_0(\overline S)-\overline S\otimes1-1\otimes\overline S
       \in\rU(\Ltw)\otimes\rU(\g[u]).
\]
We apply Lemma~\ref{lem:relative-primitives} to $\Ltw\subset\g[u]$
with the complement
\[
 \mathfrak m=\bigoplus_{r\geqslant0}\g^{(-1)^{r+1}}u^r .
\]
Then we obtain a unique decomposition $\overline S=v+w$ with
$v\in\rU(\Ltw)$ and $w\in\mathfrak m$. Since the action of
$\mathfrak t$, the involution $\omega$ and the current degree preserve
$\Ltw$ and $\mathfrak m$, the uniqueness implies that $w$ has degree
$m$, weight $\mu$ and $\omega$-eigenvalue $\epsilon$. Thus $w=xu^m$
with $x\in(\g^\epsilon)_\mu=0$, and hence
$\overline S\in\rU(\Ltw)_m$.

As $\cl$ maps the image of $\B_m$ onto $\rU(\Ltw)_m$, there exists
$Q'\in\B_m$ such that $\overline{Q'}=\overline S$. Let $Q$ be the
component of $Q'$ of weight $\mu$ and $\omega$-eigenvalue $\epsilon$.
By the decomposition of $\B$ noted above, $Q\in\B_m$. Moreover,
$\overline Q=\overline S$, because $\overline S$ has this weight and
eigenvalue. Then
\[
 S-Q=\hbar S_1,\qquad S_1\in\Y_{m-1}.
\]
If $m=0$, then $S=Q\in\B_0$, as the grading is nonnegative. If
$m>0$, then $\nabla(Q)\in\mathcal M_n\otimes\Y$ by
\eqref{eq:low-degree-module-coaction}, since
$Q\in\B_m\subset\mathcal M_n$ and $1\in\mathcal M_n$. Hence
$\hbar\nabla(S_1)\in\mathcal M_n\otimes\Y$. Because
$\mathcal M_n\otimes\Y$ is saturated, we have
$\nabla(S_1)\in\mathcal M_n\otimes\Y$. Clearly $S_1$ has the same
weight and eigenvalue as $S$. By induction, $S_1\in\B_{m-1}$, and
thus $S=Q+\hbar S_1\in\B_m$.
\end{proof}

\begin{proposition}[Saturation criterion]
\label{prop:saturation-criterion}
In the setting above, the following hold.
\begin{enumerate}[label=\textup{(\alph*)}]
\item $\B$ is $\hbar$-saturated in $\Y$, and $\cl$ induces an
isomorphism $\cli:\B/\hbar\B\xrightarrow{\sim}\rU(\Ltw)$.
\item For any ordered homogeneous basis of $\Ltw$ and any homogeneous
lifts in $\B$ of its elements, the ordered monomials in these lifts,
including $1$, form a $\C[\hbar]$-basis of $\B$. In particular, $\B$
is a homogeneous coideal quantization of $\Ltw$.
\item $\mathscr R_\nu\in\hbar^2\B_{d_\nu-2}$ for every $\nu$.
\item Let $s_\nu\in\mathscr F$ be homogeneous of degree $d_\nu-1$ with
$\pi(s_\nu)=\hbar^{-1}\mathscr R_\nu$. Then $\ker\pi$ is the
two-sided ideal of $\mathscr F$ generated by the elements
$p_\nu=r_\nu-\hbar s_\nu$. The $s_\nu$ may be chosen of weight
$\mu_\nu$, and then each $p_\nu$ is a weight vector of weight
$\mu_\nu$.
\end{enumerate}
\end{proposition}

\begin{proof}
(a) Suppose that saturation fails, and let $n$ be the least degree
with $\B_n\cap\hbar\Y_{n-1}\ne\hbar \B_{n-1}$. Choose
$z\in \B_n\cap\hbar\Y_{n-1}$ outside $\hbar \B_{n-1}$, and a
homogeneous preimage $w\in\mathscr F_n$ of $z$ under $\pi$. Since
$\overline z=0$, the reduction of $w$ modulo $\hbar$ lies in the
kernel of $\rU(\mathfrak F)\to\rU(\Ltw)$. By
Lemma~\ref{lem:classical-relations}, this reduction is a sum of
homogeneous two-sided multiples of the $r_\nu$. Lifting the factors
to $\mathscr F$ and evaluating, we obtain
\[
 z=\hbar z'+\sum_k \mathscr A_k \mathscr R_{\nu_k} \mathscr C_k,
 \qquad z'\in \B_{n-1},\qquad
 \deg \mathscr A_k+d_{\nu_k}+\deg \mathscr C_k=n,
\]
where $\mathscr A_k,\mathscr C_k\in\B$ are homogeneous and repetitions
of relations are allowed. If $d_{\nu_k}<n$, then
$\overline{\mathscr R_{\nu_k}}=0$, and by the minimality of $n$, we
have $\mathscr R_{\nu_k}\in\hbar\B_{d_{\nu_k}-1}$. If all relations of
degree $n$ in the sum also had this property, then $z$ would lie in
$\hbar\B_{n-1}$. Hence there is a relation of degree $d_\nu=n$ such
that
\begin{equation}\label{eq:minimal-homogeneous-obstruction}
\mathscr R:=\mathscr R_\nu\in \B_n\cap\hbar\Y_{n-1},\qquad
\mathscr R\notin\hbar \B_{n-1}.
\end{equation}
Put $\mu=\mu_\nu$ and $\epsilon=(-1)^n$. Then $\mathscr R$ has
weight $\mu$ and $\omega$-eigenvalue $\epsilon$, and
$(\g^\epsilon)_\mu=0$ by
Lemma~\ref{lem:classical-relation-weights}.

Since all degrees below $n$ are saturated,
Lemma~\ref{lem:low-degree-saturation} applies to $\mathcal M_n$. We
claim that
\begin{equation}\label{eq:relation-lower-coaction-pbw}
 \nabla(\mathscr R)\in\mathcal M_n\otimes\Y.
\end{equation}
We expand the coproduct of the Lie word $\mathscr R$, using
$\Delta(x)=x\otimes1+1\otimes x$ for the letters $x\in\kalg$ and
(Q3) for the letters $\sfb(v)$. Each letter $\sfb(v)$ contributes
either $\sfb(v)\otimes1$ or a tensor whose first factors lie in
$\rU(\kalg)[\hbar]$. Consider a term in which every letter $\sfb(v)$
contributes $\sfb(v)\otimes1$. If every letter $x\in\kalg$
contributes $x\otimes1$, then the term is $\mathscr R\otimes1$.
Otherwise, the term is a Lie word in elements of the form $X\otimes1$
and $1\otimes X'$, involving both forms. Such a term vanishes. Indeed,
at the lowest bracket whose two arguments involve both forms, a Lie
word in elements $X\otimes1$ is bracketed with a Lie word in elements
$1\otimes X'$, and these two commute. Thus $\mathscr R\otimes1$ is
the only nonzero term in which every letter $\sfb(v)$ contributes to
the first factor, and every remaining term has fewer than $n/e$
letters $\sfb(v)$ in its first factor. As the tensor product is
over $\C[\hbar]$, we may move powers of $\hbar$ to the second factor.
Then the first factors are products of elements of $\B$ of total
degree at most $n-e<n$. Since $1\in\mathcal M_n$, this proves
\eqref{eq:relation-lower-coaction-pbw}.

By \eqref{eq:minimal-homogeneous-obstruction}, we have
$\mathscr R=\hbar\mathscr S$ for a unique $\mathscr S\in\Y_{n-1}$.
Because $\nabla$ is linear and $\mathcal M_n\otimes\Y$ is saturated, we
obtain $\nabla(\mathscr S)\in\mathcal M_n\otimes\Y$. Moreover,
$\mathscr S$ has weight $\mu$ and $\omega$-eigenvalue $\epsilon$,
since $\hbar$ has weight zero and is fixed by $\omega$. Applying
Lemma~\ref{lem:excluded-weight-descent} with $m=n-1$, we obtain
$\mathscr S\in\B_{n-1}$, which contradicts
\eqref{eq:minimal-homogeneous-obstruction}. This proves the
saturation. As $\cl$ maps the image of $\B$ onto $\rU(\Ltw)$, we
obtain the isomorphism $\cli$.

(b) Choose an ordered homogeneous basis of $\Ltw$ and homogeneous lifts
of its elements in $\B$. Let $\{m_\lambda\}$ be the ordered monomials
in these lifts, including $1$, and put $d_\lambda=\deg m_\lambda$.
By (a), their reductions form the PBW basis of $\rU(\Ltw)$, so the
graded map
\[
 \bigoplus_\lambda\C[\hbar]e_\lambda\longrightarrow \B,
 \qquad e_\lambda\longmapsto m_\lambda,\qquad
 \deg e_\lambda=d_\lambda,
\]
is an isomorphism modulo $\hbar$. Since $\B\subset\Y$, the target is
torsion-free, and hence the map is an isomorphism by
Lemma~\ref{lem:graded-lifting}. Together with the coideal property,
this shows that $\B$ satisfies Definition~\ref{def:polynomial-coq}.

(c) By (a), $\mathscr R_\nu=\hbar\mathscr S_\nu$ with
$\mathscr S_\nu\in\B_{d_\nu-1}$, and $\mathscr S_\nu$ has
$\omega$-eigenvalue $(-1)^{d_\nu}$. On the other hand, $\omega$ acts
on $\rU(\Ltw)_{j}$ by $(-1)^{j}$, because
$\Ltw_j\subset\g^{(-1)^j}u^j$ and $\omega$ acts coefficient-wise. In particular, $\omega$ acts on $\rU(\Ltw)_{d_\nu-1}$ by
$(-1)^{d_\nu-1}$. Thus the element
$\overline{\mathscr S_\nu}\in\rU(\Ltw)_{d_\nu-1}$ vanishes. Applying
(a) again, we obtain $\mathscr S_\nu\in\hbar\B_{d_\nu-2}$.

(d) By (a), such $s_\nu$ exist, as $\pi$ is graded and surjective.
Put $Q=\mathscr F/\langle p_\nu\rangle$. Since each $p_\nu$ lies in
$\ker\pi$, the map $\pi$ factors through a graded surjection
$Q\to\B$. Because $p_\nu\equiv r_\nu$ modulo $\hbar$, it follows from
Lemma~\ref{lem:classical-relations} and (a) that
\[
 Q/\hbar Q
 \cong\rU(\mathfrak F)/\langle r_\nu\rangle
 \cong\rU(\Ltw)
 \cong \B/\hbar \B,
\]
and the composite is the map induced by $Q\to\B$, as both send
$x\mapsto x$ and $\underline v\mapsto u^ev$. Since $\B$ is
torsion-free, Lemma~\ref{lem:graded-lifting} shows that $Q\to\B$ is
an isomorphism, that is, $\ker\pi=\langle p_\nu\rangle$. Finally,
 $\mathscr F$ is the direct sum of its $\mathfrak t$-weight spaces
and $\pi$ is $\mathfrak t$-equivariant, so replacing $s_\nu$ by its
component of weight $\mu_\nu$ does not change $\pi(s_\nu)$.
\end{proof}

\section{Classification of embedded twisted Yangians}\label{sec:rigidity}

In this section and the next one, we assume that
$\omega\ne\mathrm{id}$. We show that every homogeneous coideal
quantization is generated by $\rU(\kalg)$ and a canonically normalized
lift of $u\palg$, which is determined by a pair
$(a,\K)\in\C\times\mathfrak z(\kalg)$.

\subsection{The isotropy module}

The $\kalg$-module $\palg$, with action $x\cdot y=[x,y]$, is called the
\emph{isotropy module}. The possible corrections to a normalized
degree-one lift are described by the $\kalg$-equivariant linear maps
from $\palg$ to $\g$.

\begin{lemma}\label{lem:modules}
For every nontrivial complex symmetric pair $(\g,\kalg)$
with $\g$ simple, one has
\begin{equation}\label{eq:parameter-classification}
 \Hom_{\kalg}(\palg,\kalg)=0,\qquad
 \Hom_{\kalg}(\palg,\g)
 =\C\,\mathrm{id}_{\palg}
       \oplus\ad\mathfrak z(\kalg)|_{\palg},
\end{equation}
and the map $\K\mapsto\ad(\K)|_{\palg}$ is injective on
$\mathfrak z(\kalg)$. Moreover, $\mathfrak z(\kalg)\ne0$ if and only
if $\palg$ is reducible; in this case $\dim\mathfrak z(\kalg)=1$, and
$\palg=\palg^+\oplus\palg^-$ is the sum of two nonisomorphic
irreducible $\kalg$-modules on which a suitable generator $z_0$ of
$\mathfrak z(\kalg)$ acts by $+1$ and $-1$, respectively.
\end{lemma}

\begin{proof}
For related results, see \cite[Thms.~2.2 and~2.4]{Laq84}
on invariant principal connections and \cite[\S2.3]{BDE99}
on the module structure of simple Lie triple systems.
We give a direct proof.

\emph{Step 1: preliminaries.}
Choose an $\omega$-stable compact real form $\mathfrak u$ of $\g$;
see \cite[Chap.~X, Prop.~1.4 and its proof]{Hel78}. Then $\kalg$ is
the complexification of $\mathfrak u^\omega$, so it is reductive and
its action on $\g$ is completely reducible. Indeed, the negative
Killing form on $\mathfrak u$ extends to a positive definite Hermitian
form on $\g$ for which every $\ad(x)$, $x\in\mathfrak u^\omega$, is
skew-adjoint, so the orthogonal complement of a $\kalg$-submodule is
again a $\kalg$-submodule.

The invariant bilinear form restricts nondegenerately to $\kalg$ and
$\palg$, and these subspaces are orthogonal. The action of $\kalg$ on
$\palg$ is faithful. Indeed, if $x\in\kalg$ centralizes $\palg$, then
by the Jacobi identity and $[\palg,\palg]=\kalg$, it also centralizes
$\kalg$, and hence $x\in\mathfrak z(\g)=0$. This proves the
injectivity of $\K\mapsto\ad(\K)|_\palg$.

Suppose that $0\ne V\subseteq\palg$ is a $\kalg$-submodule on which
the invariant form is nondegenerate. We claim that $V=\palg$. Put
$W=V^\perp\cap\palg$. For $v\in V$, $w\in W$ and $x\in\kalg$,
invariance gives $([v,w],x)=(v,[w,x])=0$, so $[V,W]=0$. Since
$[V,V]\subseteq\kalg$, the subspace $V+[V,V]$ is stable under
$\ad V$. It is also stable under $\kalg$ and commutes with $W$. As
$\palg=V\oplus W$, it is an ideal of $\g$, and therefore $V+[V,V]=\g$ and
$V=\palg$. Consequently, the invariant form vanishes identically on
every proper irreducible $\kalg$-submodule of $\palg$. Indeed, the
radical of its restriction is a submodule, so the restriction is
either nondegenerate or zero.

\emph{Step 2: $\palg$ reducible.}
Choose an irreducible submodule $V\subset\palg$ and, by complete
reducibility, write $\palg=V\oplus\bigoplus_j W_j$ with each $W_j$
irreducible. Since $V$ is isotropic and the form on $\palg$ is
nondegenerate, we have $(V,W_j)\ne0$ for some $j$, and we put
$W=W_j$. Then $V\cap W=0$, and the pairing $V\times W\to\C$ is
nondegenerate, because its radicals are submodules. As $V$ and $W$
are isotropic, the form is nondegenerate on the submodule
$V\oplus W$. By Step~1, we have
\[
 \palg=V\oplus W,\qquad [V,V]=[W,W]=0,
\]
where the bracket identities follow from
$([V,V],\kalg)=(V,[V,\kalg])=0$ and similarly for $W$.

Define $\partial:\g\to\g$ by $\partial|_{\kalg}=0$, $\partial|_V=1$
and $\partial|_W=-1$. Then $\partial$ is a derivation. Indeed, for
$x\in\kalg$, $v\in V$ and $w\in W$, we have
\[
 \partial[x,v]=[x,v]=[\partial x,v]+[x,\partial v],\qquad
 \partial[x,w]=-[x,w]=[\partial x,w]+[x,\partial w],
\]
both sides of the derivation identity vanish on $\kalg\times\kalg$,
$V\times V$ and $W\times W$, and because $[V,W]\subseteq\kalg$, we have
$\partial[v,w]=0=[v,w]-[v,w]=[\partial v,w]+[v,\partial w]$. Since
every derivation of a semisimple Lie algebra is inner
\cite[Chap.~II, Prop.~6.4]{Hel78}, we have $\partial=\ad(z_0)$ for
some $z_0\in\g$. As $\partial$ commutes with $\omega$, we obtain
$$\ad(\omega(z_0))=\omega \partial\omega^{-1}=\partial=\ad(z_0),$$
and hence $z_0\in\kalg$. Moreover, $[z_0,\kalg]=\partial(\kalg)=0$
implies $z_0\in\mathfrak z(\kalg)$.

Let $z\in\mathfrak z(\kalg)$. By Schur's lemma, $z$ acts on $V$ and
$W$ by some scalars $\alpha$ and $\beta$, respectively. By the
invariance of the pairing, we have
$0=([z,v],w)+(v,[z,w])=(\alpha+\beta)(v,w)$, so
$\beta=-\alpha$. Since the action is faithful, $z=\alpha z_0$. Thus
\[
 \mathfrak z(\kalg)=\C z_0,\qquad
 \palg=\palg^+\oplus\palg^-,
 \qquad \palg^+=V,\quad\palg^-=W,
\]
where $z_0$ acts by $1$ and $-1$, respectively. In particular, the
two summands are nonisomorphic. This is the Hermitian case.

Let $F:\palg\to\kalg$ be a $\kalg$-module map. For
$y_\pm\in\palg^\pm$, we have
$\pm F(y_\pm)=F([z_0,y_\pm])=[z_0,F(y_\pm)]=0$, since
$z_0\in\mathfrak z(\kalg)$. Hence $F=0$. Let $T:\palg\to\palg$ be a
$\kalg$-module map. Then $[z_0,T(y_\pm)]=\pm T(y_\pm)$, so $T$
preserves the eigenspaces $\palg^\pm$ of $\ad(z_0)|_\palg$. By
Schur's lemma, there exist $c_\pm\in\C$ such that
$T|_{\palg^\pm}=c_\pm\,\mathrm{id}_{\palg^\pm}$. Consequently,
\[
 \Hom_{\kalg}(\palg,\kalg)=0,\qquad
 T=\tfrac12(c_++c_-)\,\mathrm{id}_{\palg}
       +\tfrac12(c_+-c_-)\,\ad(z_0)|_{\palg},
\]
with uniquely determined coefficients.

\emph{Step 3: $\palg$ irreducible.}
Let $z\in\mathfrak z(\kalg)$. By Schur's lemma, $[z,y]=\alpha y$ for
some $\alpha\in\C$ and all $y\in\palg$. By invariance, we have
$0=([z,v],w)+(v,[z,w])=2\alpha(v,w)$ for $v,w\in\palg$. Since the
form on $\palg\ne0$ is nondegenerate, we have $\alpha=0$, and therefore
$z=0$ by faithfulness. Thus $\mathfrak z(\kalg)=0$. By Schur's lemma,
we also have $\End_{\kalg}(\palg)=\C\,\mathrm{id}_{\palg}$.

It remains to show that $\Hom_{\kalg}(\palg,\kalg)=0$. Write
$\kalg=\bigoplus_j\kalg_j$ as a sum of simple ideals, and suppose that
$F:\palg\to\kalg$ is a nonzero $\kalg$-module map. The projections
$\pr_j:\kalg\to\kalg_j$ are $\kalg$-module maps, so we can choose $j$
such that $F_j=\pr_j\circ F\ne0$. Because both $\palg$ and
the adjoint $\kalg$-module $\kalg_j$ are irreducible, $F_j$ is an
isomorphism. For $\ell\ne j$, $x\in\kalg_\ell$ and $y\in\palg$, we
have $F_j([x,y])=[x,F_j(y)]=0$. Hence $[\kalg_\ell,\palg]=0$, and
$\kalg_\ell=0$ by faithfulness. Therefore, $\kalg=\kalg_j$ is simple,
and $F_j^{-1}$ is a $\kalg$-module isomorphism $T:\kalg\to\palg$.
As $\kalg$ is simple, the invariant symmetric forms $(Tx,Ty)$ and
$(x,y)$ on $\kalg$ are proportional, that is, $(Tx,Ty)=c(x,y)$ for
some $c\ne0$. For $x,y,z\in\kalg$, we have
\[
 ([Tx,Ty],z)=(Tx,[Ty,z])
       =(Tx,T[y,z])=c([x,y],z),
\]
so $[Tx,Ty]=c[x,y]$. Rescaling $T$, we may assume that $c=1$. Then
the subspaces $\mathfrak a_\pm=\{x\pm Tx:x\in\kalg\}$ satisfy
\[
 [x+Tx,y-Ty]=0,\qquad
 [x\pm Tx,y\pm Ty]=2\bigl([x,y]\pm T[x,y]\bigr),
\]
so $\g=\mathfrak a_+\oplus\mathfrak a_-$ is a sum of two nonzero
commuting ideals, which is a contradiction.

\emph{Step 4: conclusion.}
Since $\g=\kalg\oplus\palg$ as $\kalg$-modules, it follows from
Steps~2 and~3 that
\[
 \Hom_{\kalg}(\palg,\g)
 =\Hom_{\kalg}(\palg,\kalg)\oplus\End_{\kalg}(\palg)
 =\C\,\mathrm{id}_{\palg}
       \oplus\ad\mathfrak z(\kalg)|_{\palg}.\qedhere
\]
\end{proof}

\subsection{Normalized lifts}

Let $B\subset\Y$ be a homogeneous coideal quantization of $\Ltw$, with
coaction $\rho=\Delta|_B$. A \emph{homogeneous lift} of $yu$,
$y\in\palg$, is an element $b(y)\in B_1$ with $\cli(b(y)\bmod\hbar
B)=yu$. Clearly lifts can be chosen linearly in $y$, by lifting a
basis of $u\palg$ and extending linearly. Since the classical element
$yu$ is primitive and $\rho$ reduces to the classical coaction, the
element $\rho(b(y))-b(y)\otimes1-1\otimes b(y)$ is divisible by
$\hbar$ in $B\otimes\Y$. The reduction of the quotient is a
first-order correction to the classical coaction. In the following
lemma, we show that this correction can be normalized to the tensor
$\gamma(yu)$ in \eqref{eq:mu-1-co}, and then the lift is determined
by equivariance.

\begin{lemma}\label{lem:lift}
Let $\omega\ne\mathrm{id}$ and let $B\subset\Y$ be a homogeneous
coideal quantization of $\Ltw$. There is a unique linear map
$\mathcal B:\palg\to B_1$ with $\overline{\mathcal B(y)}=yu$ such that
\begin{align}
 [x,\mathcal B(y)]&=\mathcal B([x,y]),\label{eq:k-equiv}\\
 \rho(\mathcal B(y))
 &=\mathcal B(y)\otimes1+1\otimes\mathcal B(y)
 +\hbar\gamma(yu)\label{eq:exact-coaction}
\end{align}
for $x\in\kalg$ and $y\in\palg$, where
$\gamma(yu)=-[1\otimes y,\Omega_{\kalg}]$. Moreover, there is
$z\in\End_{\kalg}(\palg)$ such that
\begin{equation}\label{eq:embedding-rigid}
 \mathcal B(y)
 =\mathcal J(y)-\tfrac14\hbar[y,C_{\kalg}]+\hbar z(y)
 \qquad(y\in\palg).
\end{equation}
We call $\mathcal B$ the \emph{normalized lift} of $B$.
\end{lemma}

\begin{proof}
\emph{Step 1: the raw coaction.}
Let $b:\palg\to B_1$ be a linear lift. Since $b(y)-\mathcal J(y)\in\Y_1$
has zero reduction, there is a unique linear map
$\psi^\circ:\palg\to\Y_0=\rU(\g)$ with
\[
 b(y)=\mathcal J(y)+\hbar \psi^\circ(y).
\]
Because $B$ and $\Y$ are free, the reduction of $B\otimes\Y$ is
$(B/\hbar B)\otimes_{\C}(\Y/\hbar\Y)$. The tensor
$\rho(b(y))-b(y)\otimes1-1\otimes b(y)$ has zero reduction, as
$yu$ is primitive, and hence it is divisible by $\hbar$ in
$B\otimes\Y$. By torsion-freeness, the quotient is unique. As it
has total degree zero, it lies in
$B_0\otimes\Y_0=\rU(\kalg)\otimes\rU(\g)$, that is,
\begin{equation}\label{eq:raw-exact-coaction}
 \rho(b(y))=b(y)\otimes1+1\otimes b(y)
                   +\hbar r(y),\qquad r(y)\in\rU(\kalg)\otimes\rU(\g).
\end{equation}
Comparing with \eqref{eq:coproduct-J} in $\Y\otimes\Y$ and cancelling
$\hbar$, we obtain
\begin{equation}\label{eq:reduced-F}
\nabla(\psi^\circ(y))
 =r(y)-\tfrac12[y\otimes1,\Omega_{\g}].
\end{equation}

\emph{Step 2: the Casimir correction.}
Write $\Omega_{\palg}=\Omega_{\g}-\Omega_{\kalg}$ for the tensor of the
restricted form on $\palg$. The $\kalg\otimes\palg$ and
$\palg\otimes\kalg$ components of the invariance
$[y\otimes1+1\otimes y,\Omega_\g]=0$ are
\begin{equation}\label{eq:mixed-Casimir-identities}
 [y\otimes1,\Omega_{\palg}]+[1\otimes y,\Omega_{\kalg}]=0,
 \qquad
 [y\otimes1,\Omega_{\kalg}]+[1\otimes y,\Omega_{\palg}]=0.
\end{equation}
Since $\Delta(C_{\kalg})=C_{\kalg}\otimes1+1\otimes C_{\kalg}
+2\Omega_{\kalg}$ and $y$ is primitive,
$\nabla([y,C_{\kalg}])=2[y\otimes1+1\otimes y,\Omega_{\kalg}]$. Using
\eqref{eq:mixed-Casimir-identities}, we obtain
\begin{equation}\label{eq:reduced-F-symmetric}
\nabla\!\big(-\tfrac14[y,C_{\kalg}]\big)
 =-\tfrac12\bigl([y\otimes1,\Omega_{\kalg}]
                    +[1\otimes y,\Omega_{\kalg}]\bigr)
 =\gamma(yu)-\tfrac12[y\otimes1,\Omega_{\g}].
\end{equation}

\emph{Step 3: normalization.}
Set $\mathcal T(y)=\psi^\circ(y)+\tfrac14[y,C_{\kalg}]$. Then by
\eqref{eq:reduced-F} and \eqref{eq:reduced-F-symmetric}, we have
\[
 \nabla(\mathcal T(y))=r(y)-\gamma(yu)
 \in \rU(\kalg)\otimes \rU(\g).
\]
Applying Lemma~\ref{lem:relative-primitives} to
$(\mathfrak a,\mathfrak l,\mathfrak m)=(\g,\kalg,\palg)$, we obtain
unique linear maps $\psi:\palg\to\rU(\kalg)$ and $z:\palg\to\palg$
such that $\mathcal T=\psi+z$. Put $\mathcal B(y)=b(y)-\hbar\psi(y)$.
Then $\mathcal B(y)\in B_1$, as $\rU(\kalg)=B_0$, and
$\mathcal B$ satisfies \eqref{eq:embedding-rigid}. Since $z(y)$ is
primitive, we have $\nabla(\psi(y))=r(y)-\gamma(yu)$. Subtracting
$\hbar\Delta(\psi(y))$ from \eqref{eq:raw-exact-coaction}, we obtain
\eqref{eq:exact-coaction}.

\emph{Step 4: equivariance.}
For $x\in\kalg$ and $y\in\palg$, the element
$E(x,y)=[x,\mathcal B(y)]-\mathcal B([x,y])\in B_1$ has zero
reduction. By the saturation \eqref{eq:coideal-saturation},
$E(x,y)=\hbar e(x,y)$ with $e(x,y)\in B_0=\rU(\kalg)$. On the other
hand, by \eqref{eq:embedding-rigid}, the $\kalg$-equivariance of
$\mathcal J$ and $[x,C_\kalg]=0$, we have
\[
 E(x,y)=\hbar\bigl([x,z(y)]-z([x,y])\bigr)\in\hbar\palg .
\]
Cancelling $\hbar$, we get $e(x,y)\in\rU(\kalg)\cap\palg=0$. Thus
$E=0$, which proves \eqref{eq:k-equiv} and $z\in\End_{\kalg}(\palg)$.

\emph{Step 5: uniqueness.}
Suppose that both $\mathcal B_1$ and $\mathcal B_2$ satisfy the
conditions. Then $\mathcal B_1(y)-\mathcal B_2(y)\in B_1$ has zero
reduction, and hence by saturation, it equals $\hbar\Psi(y)$ for some
$\Psi(y)\in B_0=\rU(\kalg)$. Subtracting the two coaction formulas and
cancelling $\hbar$, we see that $\Psi(y)$ is primitive, and thus
$\Psi(y)\in\kalg$ by Corollary~\ref{cor:primitives}. By
\eqref{eq:k-equiv}, $\Psi$ is $\kalg$-equivariant. Therefore,
$\Psi\in\Hom_{\kalg}(\palg,\kalg)=0$ by Lemma~\ref{lem:modules}.
\end{proof}

\begin{remark}[Semiclassical normalization]\label{rem:semiclassical}
For a linear lift $b$ of $y\mapsto yu$, the first-order coaction
tensor
\begin{equation}\label{eq:polynomial-semiclassical}
 (\cli\otimes\cl)
 \Big(
 \frac{\rho(b(y))-b(y)\otimes1-1\otimes b(y)}{\hbar}
 \bmod\hbar
 \Big)
\end{equation}
changes by $\nabla(\psi(y))$ when $b(y)$ is replaced by
$b(y)+\hbar\psi(y)$ with $\psi(y)\in\rU(\kalg)$. The proof of
Lemma~\ref{lem:lift} shows that this tensor can be normalized to
$\gamma(yu)$. Moreover, this first-order condition is equivalent to
the exact formula \eqref{eq:exact-coaction}, since the quotient in
\eqref{eq:polynomial-semiclassical} has total degree zero and
$\hbar(B\otimes\Y)$ has no component of total degree zero. As the
primitive elements of $\rU(\kalg)$ are the elements of $\kalg$, two
such normalizations differ by $\hbar\psi$ with
$\psi\in\Hom_\C(\palg,\kalg)$, and this ambiguity is removed by the
equivariance \eqref{eq:k-equiv}. In particular, in the right coideal
convention, the normalized lift of every homogeneous coideal
quantization satisfies the first-order condition of
\cite[Def.~5.3]{BR17}. Lemma~\ref{lem:lift} can be viewed as a coideal
analogue of the normalization of degree-one generators in
\cite[Thm.~3.3]{GWX25}.
\end{remark}

\subsection{Embedded realizations}

\begin{theorem}\label{thm:embedded-classification}
Let $\omega\ne\mathrm{id}$ be an involution of the complex simple Lie
algebra $\g$.  Every homogeneous coideal quantization $B\subset\Y$
of $\Ltw$ has the form
\[
 B=\C[\hbar]\left\langle \rU(\kalg),\,
 \mathcal J(y)-\tfrac14\hbar[y,C_{\kalg}]
       +\hbar\bigl(ay+[\K,y]\bigr)
 \ \middle|\ y\in\palg\right\rangle=\B_{a,\K}
\]
for unique $a\in\C$ and $\K\in\mathfrak z(\kalg)$, and its normalized
lift is $b_{a,\K}$, defined in \eqref{eq:introduction-generator}. In
particular, for fixed $(a,\K)$ there is at most one homogeneous
coideal quantization of $\Ltw$ in $\Y$ containing $b_{a,\K}(\palg)$.
\end{theorem}

\begin{proof}
Let $\mathcal B$ be the normalized lift in Lemma~\ref{lem:lift}, and
let $z\in\End_\kalg(\palg)$ be as in \eqref{eq:embedding-rigid}.
Since $\Y$ is torsion-free, $z$ is determined by $\mathcal B$. We
first show that $\rU(\kalg)$ and $\mathcal B(\palg)$ generate $B$,
and then express $z$ in terms of $(a,\K)$.

Let $B^\circ\subseteq B$ be the graded $\C[\hbar]$-subalgebra
generated by $\rU(\kalg)$ and $\mathcal B(\palg)$. The reductions of
these generators are $\kalg$ and $u\palg$, which generate $\Ltw$ by
Lemma~\ref{lem:bracket-generation}, and hence they generate
$\rU(\Ltw)$. Let $b\in B_n$ be homogeneous. We express
$\cli(\overline b)$ as a linear combination of homogeneous words in
these classical generators, and replace each letter by its lift. Then
we obtain $b_0\in B^\circ_n$ such that $b-b_0=\hbar b_1$ for some
$b_1\in B_{n-1}$. Repeating this construction, we obtain
\[
 b=b_0+\hbar b_1'+\hbar^2 b_2'+\cdots+\hbar^n b_n',
 \qquad b_r'\in B^\circ_{n-r},
\]
where the process stops because $B$ has no negative degrees. Hence
$B=B^\circ$.

By Lemma~\ref{lem:modules}, there is a unique decomposition
\begin{equation}\label{eq:z-parameterization}
 z=z_{a,\K}:=a\,\mathrm{id}_{\palg}
       +\ad(\K)|_{\palg}
 \qquad(a\in\C,\ \K\in \mathfrak z(\kalg)),
\end{equation}
and substituting it into \eqref{eq:embedding-rigid}, we obtain
$\mathcal B=b_{a,\K}$ and $B=\B_{a,\K}$.

Conversely, for every $(a,\K)$, the elements $b_{a,\K}(y)$ satisfy
\eqref{eq:k-equiv} and \eqref{eq:exact-coaction}. Indeed, $z_{a,\K}$
is $\kalg$-equivariant since $\K$ is central in $\kalg$, and
$\mathcal J$ and $y\mapsto[y,C_\kalg]$ are $\kalg$-equivariant.
Moreover, by \eqref{eq:coproduct-J} and
\eqref{eq:reduced-F-symmetric}, we have
\[
 \Delta\bigl(\mathcal J(y)-\tfrac14\hbar[y,C_{\kalg}]\bigr)
 =\bigl(\mathcal J(y)-\tfrac14\hbar[y,C_{\kalg}]\bigr)\otimes1
 +1\otimes\bigl(\mathcal J(y)-\tfrac14\hbar[y,C_{\kalg}]\bigr)
 +\hbar\gamma(yu),
\]
and $ay+[\K,y]\in\palg$ is primitive. Hence, if $b_{a,\K}(\palg)$ lies
in a homogeneous coideal quantization $B$, then by the uniqueness in
Lemma~\ref{lem:lift}, $b_{a,\K}$ is the normalized lift of $B$. Thus
$z=z_{a,\K}$, and the parameters are determined by $B$. If two
homogeneous coideal quantizations have the same parameters, then they
have the same normalized lift. Since both of them are generated by
$\rU(\kalg)$ and this lift, they coincide.

In the language of Remark~\ref{rem:realizations}, if two abstract
realizations $(B_1,\varphi_1)$ and $(B_2,\varphi_2)$ have the same
parameters, then $\varphi_1(B_1)=\varphi_2(B_2)$, and
$F=\varphi_2^{-1}\varphi_1$ is the unique embedding-preserving
isomorphism, which sends $x$ to $x$ and the normalized lift of $B_1$
to that of $B_2$.
\end{proof}

We call $a$ the spectral parameter and $\K$ the central parameter of
$B$. To obtain the uniqueness up to graded algebra isomorphism, it
remains to prove the existence for all parameters and to compare the
resulting algebras. This will be done in
Section~\ref{sec:parameter-families}.

\begin{example}\label{ex:sl2-generators}
Take $\g=\mathfrak{sl}_2$ with basis $e,f,h$ satisfying
$[h,e]=2e$, $[h,f]=-2f$, $[e,f]=h$,
and the involution $\omega(h)=h$, $\omega(e)=-e$, $\omega(f)=-f$.
Then $\kalg=\C h$ and $\palg=\C e\oplus\C f$.
The normalized form has $(e,f)=1$ and $(h,h)=2$, so
$\Omega_{\kalg}=\tfrac12 h\otimes h$ and $C_{\kalg}=\tfrac12 h^2$.
For $\K=\zeta h$, the correction in
Theorem~\ref{thm:embedded-classification} is
\[
 z(e)=(a+2\zeta)e,\qquad z(f)=(a-2\zeta)f.
\]
Denote by $B_e$ and $B_f$ the two degree-one generators in the
theorem. Then
\begin{align*}
 B_e&=\mathcal J(e)+\tfrac14\hbar(eh+he)+\hbar(a+2\zeta)e,\\
 B_f&=\mathcal J(f)-\tfrac14\hbar(fh+hf)+\hbar(a-2\zeta)f,
\end{align*}
and the normalized coactions are given by
\begin{align*}
 \Delta(B_e)&=B_e\otimes1+1\otimes B_e+\hbar h\otimes e,\\
 \Delta(B_f)&=B_f\otimes1+1\otimes B_f-\hbar h\otimes f.
\end{align*}
From these formulas, we see that $a$ is the average of the two linear
correction coefficients, and $\zeta$ is one quarter of their
difference. They also explain why the coaction alone does not
distinguish $\zeta$, since changing $\zeta$ adds primitive elements.
\end{example}

\begin{remark}\label{rem:translation}
Note that the spectral parameter depends on the reference map
$\mathcal J$. Indeed, if we replace $\mathcal J$ by
$\mathcal J+c\hbar\,\mathrm{id}_{\g}$, then the parameters of the same
embedded quantization change from $(a,\K)$ to $(a-c,\K)$. For fixed
$\mathcal J$, the spectral translation $\sigma_c$ in \eqref{eq:shift}
fixes $\kalg$ and maps $b_{a,\K}(y)$ to $b_{a+c,\K}(y)$. As
$\sigma_c$ is a graded Hopf algebra automorphism whose reduction is
the identity, it preserves the coideal, saturation and
classical-limit conditions. Thus a quantization with parameters
$(a,\K)$ exists if and only if a quantization with parameters
$(0,\K)$ exists, and in this case $\sigma_c(\B_{a,\K})=\B_{a+c,\K}$.
\end{remark}

\subsection{Left-coideal formulation}
\label{ssec:left-coideal}

The left version of Definition~\ref{def:polynomial-coq} requires
$\Delta(B^\ell)\subset\Y\otimes_{\C[\hbar]} B^\ell$, with the freeness
and classical-limit conditions unchanged. We record the
correspondence between left and right coideal quantizations in our
notation; see \cite[Eq.~(2.24)--(2.26) and \S2.3]{Lu26min} for the
Yangian anti-involution and the change of convention, and see also
\cite[\S5]{BR17} for the left convention.

Let $\vartheta$ be an involutive, form-preserving Lie algebra
anti-automorphism of $\g$ commuting with $\omega$, and let
$\mathfrak T$ be its extension to a graded algebra anti-automorphism
of $\Y$ \cite[\S2.3]{Lu26min}:
\[
 \mathfrak T(x)=\vartheta(x),\qquad
 \mathfrak T(\mathcal J(x))=\mathcal J(\vartheta x),\qquad
 \mathfrak T(\hbar)=\hbar,\qquad
 \Delta\circ\mathfrak T=(\mathfrak T\otimes\mathfrak T)\circ\Delta^{\mathrm{op}}.
\]
For example, one may take $\vartheta=-\mathrm{id}_{\g}$. In the split
and quasi-split setting of Section~\ref{sec:minimal}, we use the
compatible Chevalley anti-involution given by
$\vartheta(\xi_i)=\xi_i$ and $\vartheta(x_i^\pm)=x_i^\mp$.

\begin{proposition}\label{prop:left-coideal}
Let $B\subset\Y$ be a homogeneous right coideal quantization of
$\Ltw$. Then the map
\[
 \varphi^\ell:B^{\mathrm{op}}\hookrightarrow\Y,\qquad
 \varphi^\ell(b^{\mathrm{op}})=\mathfrak T(b),
\]
realizes $B^{\mathrm{op}}$ as a homogeneous left coideal quantization.
This correspondence is involutive and preserves embedding-preserving
isomorphisms. If $\mathcal B$ is the normalized lift of $B$ with
parameters $(a,\K)$, then
$\mathcal B^\ell(y)=\mathcal B(\vartheta y)^{\mathrm{op}}$ satisfies
\begin{align}
 \varphi^\ell(\mathcal B^\ell(y))
 &=\mathcal J(y)+\tfrac14\hbar[y,C_{\kalg}]
       +\hbar\bigl(ay+[y,\K^\ell]\bigr),
 \qquad \K^\ell=\vartheta(\K),
 \label{eq:left-embedding-rigid}\\
 \Delta(\varphi^\ell(\mathcal B^\ell(y)))
 &=\varphi^\ell(\mathcal B^\ell(y))\otimes1
       +1\otimes\varphi^\ell(\mathcal B^\ell(y))
       +\hbar[y\otimes1,\Omega_{\kalg}].
 \label{eq:left-exact-coaction}
\end{align}
\end{proposition}

\begin{proof}
Since $\mathfrak T$ is an injective graded anti-homomorphism,
$\varphi^\ell$ is an injective graded homomorphism. By the identity
$\Delta\mathfrak T=(\mathfrak T\otimes\mathfrak T)\Delta^{\mathrm{op}}$,
we have $\Delta(\mathfrak T(B))\subset\Y\otimes\mathfrak T(B)$. As
$\mathfrak T(\hbar\Y)=\hbar\Y$, the image $\mathfrak T(B)$ is
saturated and free. Modulo $\hbar$, $\mathfrak T$ induces the
anti-automorphism of $\rU(\g[u])$ determined by
$xu^r\mapsto\vartheta(x)u^r$. Since $\vartheta$ commutes with
$\omega$, it maps $\rU(\Ltw)$ onto itself. This proves the first
statement. The correspondence is involutive because
$\mathfrak T^2=\mathrm{id}$, and an embedding-preserving isomorphism
$F$ induces the embedding-preserving isomorphism $F^{\mathrm{op}}$.

To prove \eqref{eq:left-embedding-rigid}, we apply $\mathfrak T$ to
\eqref{eq:embedding-rigid} at $\vartheta y$. As $\{\vartheta v_r\}$
and $\{\vartheta v^r\}$ are dual bases of $\kalg$, we have
$\mathfrak T(C_\kalg)=\sum_r\vartheta(v^r)\vartheta(v_r)=C_\kalg$, and
hence $\mathfrak T([\vartheta y,C_\kalg])=[C_\kalg,y]=-[y,C_\kalg]$.
Moreover, $\mathfrak T([\K,\vartheta y])=[y,\vartheta\K]$. To prove
\eqref{eq:left-exact-coaction}, we apply
$(\mathfrak T\otimes\mathfrak T)\circ\operatorname{flip}$ to
\eqref{eq:exact-coaction} at $\vartheta y$. Using the symmetry of
$\Omega_\kalg$ and $(\vartheta\otimes\vartheta)\Omega_\kalg=\Omega_\kalg$,
we obtain
\[
 (\mathfrak T\otimes\mathfrak T)
 \bigl(-[\vartheta y\otimes1,\Omega_\kalg]\bigr)
 =-\sum_r[\vartheta v_r,y]\otimes\vartheta v^r
 =[y\otimes1,\Omega_\kalg].\qedhere
\]
\end{proof}

\section{Existence and parameter isomorphisms}
\label{sec:parameter-families}

By Theorem~\ref{thm:embedded-classification}, the generators of any
homogeneous coideal quantization are determined. In this section,
using the saturation criterion in Section~\ref{sec:criterion}, we
show that for every value of the parameters, the subalgebras
generated by these elements have the required classical limit, and
that their isomorphism type does not depend on the parameters.
Throughout this section, $(\g,\kalg)$ is an arbitrary nontrivial
complex symmetric pair.

\subsection{The candidate coideal subalgebras}

For $a\in\C$ and $\K\in \mathfrak z(\kalg)$, recall from
\eqref{eq:introduction-generator} the elements
$b_{a,\K}(y)=\mathcal J(y)-\tfrac14\hbar[y,C_{\kalg}]
 +a\hbar y+\hbar[\K,y]$, and put
\begin{equation}\label{eq:canonical-parameter-coideal}
 \B_{a,\K}
 :=\C[\hbar]\bigl\langle \rU(\kalg),\
 b_{a,\K}(\palg)\bigr\rangle
 \subset\Y.
\end{equation}
By the $\kalg$-equivariance of $\mathcal J$, the centrality of
$C_{\kalg}$ and $\K$ in $\rU(\kalg)$, and the Yangian coproduct
\eqref{eq:coproduct-J}, a direct calculation shows that for
$x\in\kalg$ and $y\in\palg$,
\begin{equation}\label{eq:direct-pbw-generator-properties}
 \begin{aligned}
 [x,b_{a,\K}(y)]&=b_{a,\K}([x,y]),\\
 \Delta(b_{a,\K}(y))&=b_{a,\K}(y)\otimes1+1\otimes b_{a,\K}(y)
       -\hbar[1\otimes y,\Omega_{\kalg}]\\
       &=b_{a,\K}(y)\otimes1+1\otimes b_{a,\K}(y)
       -\hbar\sum_r v_r\otimes[y,v^r],
 \end{aligned}
\end{equation}
where $v_r,v^r$ are dual bases of $\kalg$. In particular, $\B_{a,\K}$
is a graded right coideal subalgebra. Apart from
$b_{a,\K}(y)\otimes1$, every term on the right-hand side of the second
formula has its first factor in $\rU(\kalg)[\hbar]$. Finally, since
the invariant form and $C_\kalg$ are preserved by $\omega$ and $\K$ is
fixed by $\omega$, we have
\[
 \omega(b_{a,\K}(y))=-b_{a,\K}(y)\qquad(y\in\palg).
\]

\subsection{The affine datum}

For nontrivial involutions, the classical data needed in
Section~\ref{sec:criterion} are provided by the following lemma.

\begin{lemma}
\label{lem:affine-degree-one}
There is an affine datum for $(\g,\omega)$, in the sense of
Definition~\ref{def:affine-datum}, with $e=1$.
\end{lemma}

\begin{proof}
We may replace $\omega$ by any conjugate $\phi\omega\phi^{-1}$,
$\phi\in\operatorname{Aut}(\g)$. Indeed, composing an affine datum for
$\phi\omega\phi^{-1}$ with the coefficient-wise action of $\phi^{-1}$,
which maps the twisted loop algebra of $\phi\omega\phi^{-1}$ onto $L$
and its fixed-point algebra and $(-1)$-eigenspace onto $\kalg$ and
$\palg$, we obtain an affine datum for $\omega$.

By \cite[Thm.~8.6]{Kac90}, after such a conjugation $\omega$ is the
automorphism of type $(s;r)$ for a sequence of nonnegative integers
$(s_i)_{i\in\hI}$ satisfying
\[
 r\sum_i m_i s_i=2,\qquad \gcd(s_i)=1,
\]
where the $m_i$ are the positive integer marks of the affine Dynkin
diagram and $r\in\{1,2\}$ is the order of the image of $\omega$ modulo
inner automorphisms. The loop realization in
\cite[proof of Thm.~8.5]{Kac90} gives an isomorphism
$\lambda:\widehat\g'/\C\mathsf c\to L$ satisfying (A1) with these
$s_i$. If $s_i\geqslant2$ for some $i$, then the displayed sum forces
$r=m_i=1$, $s_i=2$ and $s_j=0$ for $j\ne i$, which contradicts
$\gcd(s_i)=1$. Hence each $s_i$ is $0$ or $1=e$.

Property (A2) follows from \cite[Prop.~8.6(b)]{Kac90}; in fact,
$\mathfrak t$ is a Cartan subalgebra of $\kalg$. By
\cite[Prop.~8.6(c)]{Kac90} and the
identification of $\palg$ with $u\palg=L_1$ as $\kalg$-modules,
\[
 \palg=\bigoplus_{i\in \hI^1}\rU(\kalg)\cdot y_i
\]
with each summand irreducible. For $j\in\hI^0$ and $i\in\hI^1$, the
affine relation $[f_j,e_i]=0$ implies $[\overline f_j,y_i]=0$. Thus
$y_i$ is annihilated by all negative simple root vectors of $\kalg$,
that is, $y_i$ is a lowest weight vector. By the triangular
decomposition of $\rU(\kalg)$, the module $\rU(\kalg)\cdot y_i$ is
spanned by the iterated brackets of the $\overline e_j$,
$j\in\hI^0$, applied to $y_i$. This proves (A3).
\end{proof}

In the rest of this section, we fix such an affine datum, with Serre
words $r_\nu$ of degrees
$d_\nu\geqslant2$ and $\mathfrak t$-weights $\mu_\nu$, as in
Lemmas~\ref{lem:classical-relations}
and~\ref{lem:classical-relation-weights}.

\begin{example}\label{ex:sl2-serre}
For the pair $(\mathfrak{sl}_2,\C h)$ of
Example~\ref{ex:sl2-generators}, the affine Cartan matrix is
$\begin{pmatrix}2&-2\\-2&2\end{pmatrix}$, with $s_0=s_1=1$,
$\overline e_1=eu$ and $\overline e_0=fu$. Write
$E=\underline e$ and $F=\underline f$, so that $\pi_0(E)=eu$ and
$\pi_0(F)=fu$. The index set $\mathcal S$ of
\eqref{eq:serre-index-set} is $\{(1,0),(0,1)\}$, and the corresponding
Serre words are
\[
 r_+=[E,[E,[E,F]]],\qquad r_-=[F,[F,[F,E]]].
\]
The image of $r_+$ vanishes, since
\[
 [eu,fu]=hu^2,\qquad [eu,hu^2]=-2eu^3,\qquad [eu,-2eu^3]=0.
\]
The word $r_+$ has degree four and weight $3\cdot2-2=4$ with respect
to $h$. As four is even, the relevant coefficient space in
\eqref{eq:classical-serre-weight-exclusion} is $\g^+=\C h$, which only
has weight zero, and hence $(\g^+)_4=0$. Similarly, $r_-$ has degree
four and weight $-4$. By Lemma~\ref{lem:classical-relations}, $r_\pm$
generate the full classical kernel, while the bracket calculation
above only shows that they lie in it. Note that the weight exclusion
holds for these chosen generators, but not for every element of the
kernel. For instance, $[F,r_+]$ is a relation of degree five and
weight two, and $\g^-=\palg$ contains the weight-two vector $e$.
\end{example}

\subsection{Saturation and the PBW basis}

\begin{theorem}
\label{thm:coideal-pbw}
For every nontrivial complex symmetric pair and every
$a\in\C$, $\K\in\mathfrak z(\kalg)$, the subalgebra
$\B_{a,\K}$ is $\hbar$-saturated, its natural
reduction is $\rU(\Ltw)$, and it is free over $\C[\hbar]$ with the
PBW bases of Proposition~\ref{prop:saturation-criterion}(b).
In particular, $\B_{a,\K}$ is a twisted Yangian
associated with $(\g,\kalg)$. Moreover, every Serre word satisfies
$\pi(r_\nu)\in\hbar^2\B_{a,\K}$, where $\pi$ is the evaluation
\eqref{eq:quantum-word-evaluation} with $\sfb=b_{a,\K}$.
\end{theorem}

\begin{proof}
We apply Proposition~\ref{prop:saturation-criterion} with $e=1$,
$V=\palg$, the affine datum in Lemma~\ref{lem:affine-degree-one} and
$\sfb=b_{a,\K}$. By \eqref{eq:direct-pbw-generator-properties}, the
map $\sfb$ is $\kalg$-equivariant. Condition (Q1) holds since
$\overline{b_{a,\K}(y)}=yu$, (Q2) follows from
$\omega(b_{a,\K}(y))=-b_{a,\K}(y)$, and (Q3) is the second formula of
\eqref{eq:direct-pbw-generator-properties}.
\end{proof}

\subsection{Parameter isomorphisms and uniqueness}

By Theorem~\ref{thm:coideal-pbw}, homogeneous coideal quantizations
exist for all parameters $(a,\K)$. For the same $\K$ and different
$a$, these subalgebras are related by spectral translations; see
Remark~\ref{rem:translation}. To compare two different values of
$\K$, we apply a scalar-valued algebra homomorphism to the first
factor of the coproduct, which changes $b_{a,\K}(y)$ by the required
term. For this, we first need to check that the scalar assignment
respects all relations, and here we use
Proposition~\ref{prop:saturation-criterion}(d).

\begin{lemma}
\label{lem:scalar-evaluation}
For $a\in\C$ and $\K,\mathscr D\in\mathfrak z(\kalg)$, the assignments
\begin{equation}\label{eq:general-scalar-evaluation}
 \chi_{a,\K}^{\mathscr D}(x)=(\mathscr D,x),\quad
 \chi_{a,\K}^{\mathscr D}(b_{a,\K}(y))=0,\quad
 \chi_{a,\K}^{\mathscr D}(\hbar)=\hbar\qquad (x\in\kalg,\ y\in\palg)
\end{equation}
define a graded $\C[\hbar]$-algebra homomorphism
$\chi_{a,\K}^{\mathscr D}:\B_{a,\K}\to\C[\hbar]$.
\end{lemma}

\begin{proof}
Write $\B=\B_{a,\K}$ and $\chi=\chi_{a,\K}^{\mathscr D}$. Since
$\mathscr D$ is central and the form is invariant, we have
$(\mathscr D,[x,x'])=([\mathscr D,x],x')=0$ for $x,x'\in\kalg$. As
the target is commutative, the assignments $x\mapsto(\mathscr D,x)$,
$\underline y\mapsto0$, $\hbar\mapsto\hbar$ respect the defining
relations of $\mathfrak F$ (with $V=\palg$), and hence define a graded
homomorphism $\widetilde\chi:\mathscr F\to\C[\hbar]$, where
$\mathscr F=\C[\hbar]\otimes\rU(\mathfrak F)$. Because
$\pi:\mathscr F\twoheadrightarrow\B$ is surjective,
$\widetilde\chi$ descends to $\B$ if and only if it kills
$\ker\pi$.

By Proposition~\ref{prop:saturation-criterion}(d), $\ker\pi$ is
generated by the elements $p_\nu=r_\nu-\hbar s_\nu$, where $p_\nu$ has
weight $\mu_\nu$. By Lemma~\ref{lem:classical-relation-weights}, we
have $\mu_\nu\ne0$. For $t\in\mathfrak t$ with $\mu_\nu(t)\ne0$, since
$\C[\hbar]$ is commutative, we have
\[
 \mu_\nu(t)\widetilde \chi(p_\nu)
 =\widetilde \chi([t,p_\nu])
 =[\widetilde \chi(t),\widetilde \chi(p_\nu)]=0.
\]
Hence $\widetilde\chi(p_\nu)=0$ for all $\nu$, and $\widetilde\chi$
descends to the required homomorphism $\chi$.
\end{proof}

\begin{lemma}
\label{lem:central-translation-quantization}
For $a\in\C$ and $\K,\mathscr D\in\mathfrak z(\kalg)$, the map
\[
 \varrho^{\mathscr D}_{a,\K}:=(\chi^{\mathscr D}_{a,\K}\otimes\mathrm{id})\circ\Delta:
 \B_{a,\K}\longrightarrow \B_{a,\K+\mathscr D}
\]
is a graded algebra isomorphism, with
\begin{equation}\label{eq:direct-central-translation}
 \varrho^{\mathscr D}_{a,\K}(x)=x+(\mathscr D,x)1,\qquad
 \varrho^{\mathscr D}_{a,\K}(b_{a,\K}(y))=b_{a,\K+\mathscr D}(y)
 \qquad (x\in\kalg,\ y\in\palg).
\end{equation}
Its inverse is $\varrho^{-\mathscr D}_{a,\K+\mathscr D}$, and it preserves
the right $\Y$-coaction:
$\Delta\circ\varrho^{\mathscr D}_{a,\K}
=(\varrho^{\mathscr D}_{a,\K}\otimes\mathrm{id})\circ\Delta$.
\end{lemma}

\begin{proof}
By the coideal property and Lemma~\ref{lem:scalar-evaluation},
$\varrho^{\mathscr D}_{a,\K}$ is a graded algebra homomorphism
$\B_{a,\K}\to\Y$. For dual bases $v_r,v^r$ of $\kalg$, by the identity
$\sum_r(\mathscr D,v_r)v^r=\mathscr D$ and
\eqref{eq:direct-pbw-generator-properties}, we have
\begin{align*}
 \varrho^{\mathscr D}_{a,\K}(x)&=x+(\mathscr D,x)1,\\
 \varrho^{\mathscr D}_{a,\K}(b_{a,\K}(y))
 &=b_{a,\K}(y)-\hbar\sum_r(\mathscr D,v_r)[y,v^r]
 =b_{a,\K}(y)+\hbar[\mathscr D,y]=b_{a,\K+\mathscr D}(y).
\end{align*}
These images generate $\B_{a,\K+\mathscr D}$. Similarly, the map
$\varrho^{-\mathscr D}_{a,\K+\mathscr D}$ subtracts the scalar shift
and sends $b_{a,\K+\mathscr D}(y)$ back to $b_{a,\K}(y)$. Hence both
composites fix $\hbar$ and the generators, and the two maps are
inverse to each other. Finally, by coassociativity, we have
\begin{align*}
 \Delta \varrho^{\mathscr D}_{a,\K}
 &=(\chi^{\mathscr D}_{a,\K}\otimes\mathrm{id}\otimes\mathrm{id})
       (\mathrm{id}\otimes\Delta)\Delta\\
 &=(\chi^{\mathscr D}_{a,\K}\otimes\mathrm{id}\otimes\mathrm{id})
       (\Delta\otimes\mathrm{id})\Delta
   =(\varrho^{\mathscr D}_{a,\K}\otimes\mathrm{id})\Delta.\qedhere
\end{align*}
\end{proof}

\begin{remark}\label{rem:central-vs-spectral}
Since $\varepsilon(\varrho^{\mathscr D}_{a,\K}(x))=(\mathscr D,x)$ for
$x\in\kalg$, the central translation is not the restriction of a Hopf
algebra automorphism of $\Y$ if $\mathscr D\ne0$. This is different
from the spectral translations $\sigma_c$. In the Hermitian case, the
characters $\chi^{\mathscr D}_{a,\K}$ have the same form as the counit
in \cite[Thm.~5.5]{BR17}, whose value on the center of $\kalg$ is a
free parameter there. For $\imath$quantum groups, the isomorphisms
between algebras with different parameters in \cite[Lem.~2.5.1]{Wat21}
also consist of the restriction of a Hopf algebra automorphism and a
part which does not come from the quantum group.
\end{remark}

\begin{theorem}[Uniqueness up to isomorphism]
\label{thm:main-uniqueness}
For every nontrivial complex symmetric pair $(\g,\kalg)$,
the associated twisted Yangian exists and is unique up to
graded $\C[\hbar]$-algebra isomorphism.
More explicitly, let $a_0,a_1\in\C$ and $\K_0,\K_1\in\mathfrak z(\kalg)$,
and put $\mathscr D=\K_1-\K_0$.
There is a graded $\C[\hbar]$-algebra isomorphism
\begin{equation}\label{eq:central-parameter-isomorphism}
 \Phi_{a_0,\K_0}^{a_1,\K_1}:\B_{a_0,\K_0}
       \xrightarrow{\ \sim\ }\B_{a_1,\K_1},\qquad
 x\longmapsto x+(\mathscr D,x)1,\quad b_{a_0,\K_0}(y)\longmapsto b_{a_1,\K_1}(y),
\end{equation}
for $x\in\kalg$ and $y\in\palg$, and it satisfies
\[
 \Delta\circ\Phi_{a_0,\K_0}^{a_1,\K_1}
 =\bigl(\Phi_{a_0,\K_0}^{a_1,\K_1}\otimes\sigma_{a_1-a_0}\bigr)\circ\Delta .
\]
\end{theorem}

\begin{proof}
By Theorem~\ref{thm:coideal-pbw}, there is a quantization for each
$(a,\K)$, and by Theorem~\ref{thm:embedded-classification}, these are
all the homogeneous coideal quantizations. Define
\[
 \Phi_{a_0,\K_0}^{a_1,\K_1}
 :=\sigma_{a_1-a_0}\circ\varrho^{\mathscr D}_{a_0,\K_0},
\]
that is, the composite of the central translation in
Lemma~\ref{lem:central-translation-quantization} and the spectral
translation \eqref{eq:shift}; the latter maps $\B_{a_0,\K_1}$ onto
$\B_{a_1,\K_1}$ by Remark~\ref{rem:translation}. Clearly
$\Phi_{a_0,\K_0}^{a_1,\K_1}$ is given by the assignments
\eqref{eq:central-parameter-isomorphism}. Since
$\sigma=\sigma_{a_1-a_0}$ is a Hopf algebra automorphism and
$\varrho=\varrho^{\mathscr D}_{a_0,\K_0}$ preserves the coaction, we
have
\[
 \Delta\sigma\varrho=(\sigma\otimes\sigma)\Delta\varrho
 =(\sigma\otimes\sigma)(\varrho\otimes\mathrm{id})\Delta
 =(\sigma\varrho\otimes\sigma)\Delta.\qedhere
\]
\end{proof}

\begin{remark}\label{rem:BR17-uniqueness}
Consider a homogeneous quantization in the sense of
\cite[Def.~5.3]{BR17}, realized inside the Yangian as in
\cite[Thms.~5.5 and~5.6]{BR17}. After passing to the right coideal
convention in \S\ref{ssec:left-coideal}, it is a homogeneous coideal
quantization in the sense of Definition~\ref{def:polynomial-coq}.
Hence Theorem~\ref{thm:main-uniqueness}, together with
Theorem~\ref{thm:identity-uniqueness} for the identity involution,
recovers the uniqueness statements of \cite[Thms.~5.5 and~5.6]{BR17}
in the homogeneous setting, without the condition on the first-order
coaction.
\end{remark}

\begin{corollary}\label{cor:HL26}
For $c\in\C$ and $\mathcal K\in\mathfrak z(\kalg)$, the twisted Yangian
of \cite[Def.~2.8]{HL26}, generated in our normalization by $\kalg$
and the elements
$\mathcal J(y)-\tfrac14\hbar[y,C_{\kalg}]-c\hbar[y,\mathcal K]$,
$y\in\palg$, equals $\B_{0,c\mathcal K}$. In particular, it is a
homogeneous coideal quantization of $\Ltw$. In nontrivial quasi-split
type $\mathsf A_{2n}$, the $\imath$Yangian, which is identified with
$\B_{0,\K_\circ}$ in \cite[Thm.~3.4(3)]{HL26} (see
\eqref{eq:A2n-central-element}), is isomorphic as a graded algebra to
the traditional twisted Yangian $\B_{0,0}$.
\end{corollary}

\begin{proof}
Since $-c\hbar[y,\mathcal K]=\hbar[c\mathcal K,y]$, the generators
are $b_{0,c\mathcal K}(y)$, and Theorem~\ref{thm:coideal-pbw} applies.
The last statement follows from Theorem~\ref{thm:main-uniqueness}.
\end{proof}

\section{The identity involution}\label{sec:identity-involution}

In this section, we assume that $\omega=\mathrm{id}$, so that
\[
 \kalg=\g,\qquad \palg=0,\qquad
 \Ltw=\g[u^2],
\]
and $(e,V)=(2,\g)$ in the notation of Section~\ref{sec:criterion}. We
keep the convention $\deg u=\deg\hbar=1$, so the first positive-degree
generators have degree two. We construct homogeneous coideal quantizations of
$\g[u^2]$, prove their PBW property, and classify all embedded
realizations.

\subsection{The degree-two generators}

Since $\g[u^2]$ is generated by $\g$ and $u^2\g$, we first construct
lifts of $u^2\g$ with the required coideal property. Recall the dual
bases $\{x_r\}$, $\{x^r\}$ of $\g$, $\Omega_\g=\sum_rx_r\otimes x^r$,
and put $C_{\g}=\sum_r x_rx^r$. Let $c_{\g}\ne0$ be the scalar by
which $C_{\g}$ acts on the adjoint module, so that
\begin{equation}\label{eq:identity-adjoint-casimir}
 \sum_r[x_r,[x^r,x]]=c_{\g}x
 \qquad(x\in\g).
\end{equation}
Define linear maps $\mathcal Q,\mathcal G:\g\to\Y_2$ by
\begin{equation}\label{eq:identity-degree-two-generator}
 \mathcal Q(x)=\frac1{c_{\g}}
       \sum_r[\mathcal J(x_r),\mathcal J([x^r,x])],\qquad
 \mathcal G(x)=\mathcal Q(x)
       -\tfrac14\hbar[\mathcal J(x),C_{\g}].
\end{equation}
Both maps do not depend on the choice of dual bases, and by
\eqref{eq:identity-adjoint-casimir}, we have
\[
 \overline{\mathcal Q(x)}
 =\frac1{c_{\g}}\sum_r[x_r u,[x^r,x]u]=xu^2,
 \qquad \overline{\mathcal G(x)}=xu^2.
\]
The correction term in $\mathcal G$ cancels the coproduct terms
containing $\mathcal J(\g)$ in the first factor. Note that although
$C_{\g}$ is central in $\rU(\g)$, its commutator with $\mathcal J(x)$
is taken in $\Y$ and is nonzero in general. The element
$\mathcal G(x)$ is the right coideal version of the degree-two
generator computed in \cite[\S7.4]{BR17}. Below we give the
calculation for every simple $\g$, including $\mathfrak{sl}_2$.

For $a\in\C$, put
\begin{equation}\label{eq:identity-parameter-generator}
 \mathcal G_a(x)=\mathcal G(x)+2a\hbar\mathcal J(x)
                                      +a^2\hbar^2x,
\end{equation}
and let $\B^{\mathrm{ev}}_a$ be the $\C[\hbar]$-subalgebra generated
by $\g$ and $\mathcal G_a(\g)$. As the last term belongs to
$\rU(\g)[\hbar]$, it does not affect this subalgebra. We include it in
order to obtain the exact formula for spectral translations below.

\begin{lemma}\label{lem:identity-coideal-generators}
For $a,c\in\C$ and $x,z\in\g$, one has
\begin{gather}
 [z,\mathcal G_a(x)]=\mathcal G_a([z,x]),\qquad
 \overline{\mathcal G_a(x)}=xu^2,
 \label{eq:identity-equivariance-limit}\\
 \sigma_c(\mathcal G_a(x))=\mathcal G_{a+c}(x).
 \label{eq:identity-spectral-translation}
\end{gather}
Moreover,
\begin{equation}\label{eq:identity-generator-coproduct}
 \Delta(\mathcal G_a(x))
  =\mathcal G_a(x)\otimes1+1\otimes\mathcal G_a(x)
    -\hbar[1\otimes(\mathcal J(x)+a\hbar x),\Omega_{\g}]
       +\hbar^2R(x),
\end{equation}
where $R(x)\in\rU(\g)\otimes\rU(\g)$ is independent of $a$. In
particular, $\B^{\mathrm{ev}}_a$ is a graded right coideal subalgebra,
and $\sigma_c(\B^{\mathrm{ev}}_a)=\B^{\mathrm{ev}}_{a+c}$.
\end{lemma}

\begin{proof}
By the equivariance of $\mathcal J$ and the invariance of
$\Omega_{\g}$, the map $\mathcal Q$ is $\g$-equivariant. Since
$[z,C_{\g}]=0$, it follows from the Jacobi identity that
$[z,[\mathcal J(x),C_{\g}]]=[\mathcal J([z,x]),C_{\g}]$, and hence
$\mathcal G$ and $\mathcal G_a$ are also $\g$-equivariant. Their
reductions are computed above.

The spectral translation $\sigma_c$ fixes $\g$ and sends
$\mathcal J(x)$ to $\mathcal J(x)+c\hbar x$. Expanding the commutators
in $\mathcal Q$ and using \eqref{eq:identity-adjoint-casimir} twice
(note $[\mathcal J(x_r),[x^r,x]]=\mathcal J([x_r,[x^r,x]])$ and
$[x_r,\mathcal J([x^r,x])]=\mathcal J([x_r,[x^r,x]])$), we obtain
\[
 \sigma_c(\mathcal Q(x))
   =\mathcal Q(x)+2c\hbar\mathcal J(x)+c^2\hbar^2x.
\]
The Casimir correction is unchanged, as $[x,C_{\g}]=0$. Hence
$\mathcal G_a=\sigma_a\circ\mathcal G$, which proves
\eqref{eq:identity-spectral-translation}. Since $\sigma_a$ is a Hopf
algebra automorphism, it suffices to compute the coproduct for $a=0$.

Set
\[
 \Pi(x)=\mathcal J(x)\otimes1+1\otimes\mathcal J(x),
 \qquad \Theta(x)=[x\otimes1,\Omega_{\g}],
\]
so that $\Delta(\mathcal J(x))=\Pi(x)+\tfrac12\hbar\Theta(x)$.
As $[\Pi(x),\Pi(x')]=[\mathcal J(x),\mathcal J(x')]\otimes1
+1\otimes[\mathcal J(x),\mathcal J(x')]$, expanding the commutators,
we obtain
\begin{align*}
 \nabla(\mathcal Q(x))
 ={}&\frac{\hbar}{2c_{\g}}\sum_r
   \bigl([\Pi(x_r),\Theta([x^r,x])]
            +[\Theta(x_r),\Pi([x^r,x])]\bigr)
 +\frac{\hbar^2}{4c_{\g}}
       \sum_r[\Theta(x_r),\Theta([x^r,x])].
\end{align*}
To simplify the linear term, note that by invariance, we have
$\sum_r x_r\otimes[x^r,x]=[x\otimes1,\Omega_\g]$, and its flip is
$[1\otimes x,\Omega_\g]=-[x\otimes1,\Omega_\g]$. Hence this tensor is
skew-symmetric. For a skew-symmetric tensor $\sum_ka_k\otimes b_k$, we
have
\[
 \sum_k\ad(a_k)\ad(b_k)=\tfrac12\ad\Bigl(\sum_k[a_k,b_k]\Bigr),
\]
and therefore by \eqref{eq:identity-adjoint-casimir},
\[
 \sum_r\ad(x_r)\ad([x^r,x])=\tfrac12c_{\g}\ad(x).
\]
Writing $\Theta(w)=\sum_k[w,x_k]\otimes x^k=-\sum_kx_k\otimes[w,x^k]$
and using the equivariance of $\mathcal J$, we obtain
\begin{align*}
 \sum_r[\mathcal J(x_r)\otimes1,\Theta([x^r,x])]
 &=\sum_{r,k}\mathcal J([x_r,[[x^r,x],x_k]])\otimes x^k
  =\tfrac12c_\g[\mathcal J(x)\otimes1,\Omega_\g],\\
 \sum_r[1\otimes\mathcal J(x_r),\Theta([x^r,x])]
 &=-\sum_{r,k}x_k\otimes\mathcal J([x_r,[[x^r,x],x^k]])
  =-\tfrac12c_\g[1\otimes\mathcal J(x),\Omega_\g].
\end{align*}
Thus $\sum_r[\Pi(x_r),\Theta([x^r,x])]
=\tfrac12c_{\g}[\mathcal J(x)\otimes1-1\otimes\mathcal J(x),\Omega_{\g}]$.
By the skew-symmetry of $\sum_rx_r\otimes[x^r,x]$, the other mixed sum
$\sum_r[\Theta(x_r),\Pi([x^r,x])]=-\sum_r[\Pi([x^r,x]),\Theta(x_r)]$
has the same value. Consequently,
\begin{equation}\label{eq:identity-Q-reduced-coproduct}
 \nabla(\mathcal Q(x))
  =\tfrac12\hbar
       [\mathcal J(x)\otimes1-1\otimes\mathcal J(x),
          \Omega_{\g}]
    +\frac{\hbar^2}{4c_{\g}}
       \sum_r[\Theta(x_r),\Theta([x^r,x])].
\end{equation}

On the other hand, by
$\Delta(C_{\g})=C_{\g}\otimes1+1\otimes C_{\g}+2\Omega_{\g}$, we have
\[
 \nabla([\mathcal J(x),C_{\g}])
 =2[\Pi(x),\Omega_{\g}]
       +\hbar[\Theta(x),\Omega_{\g}],
\]
since $\Theta(x)\in\g\otimes\g$ commutes with
$C_{\g}\otimes1+1\otimes C_{\g}$. Subtracting $\tfrac14\hbar$ times
this identity from \eqref{eq:identity-Q-reduced-coproduct}, we obtain
\[
 \nabla(\mathcal G(x))
   =-\hbar[1\otimes\mathcal J(x),\Omega_{\g}]
                                      +\hbar^2R(x),
 \qquad
 R(x)=\frac1{4c_{\g}}\sum_r[\Theta(x_r),\Theta([x^r,x])]
                            -\tfrac14[\Theta(x),\Omega_{\g}],
\]
and every tensor factor of $R(x)$ lies in $\rU(\g)$. Applying
$\sigma_a\otimes\sigma_a$ proves
\eqref{eq:identity-generator-coproduct}.

Finally,
$[1\otimes\mathcal J(x),\Omega_{\g}]=\sum_r x_r\otimes\mathcal J([x,x^r])$,
so all correction terms in \eqref{eq:identity-generator-coproduct} have
first factors in $\rU(\g)[\hbar]\subset\B^{\mathrm{ev}}_a$. Together
with the primitivity of $\g$, this proves the right coideal property.
\end{proof}

\subsection{Saturation and the PBW basis}

We use the untwisted affine Lie algebra $\widehat\g$ of $\g$ and its
loop realization $\widehat{\g}'/\C\mathsf c\cong\g[v,v^{-1}]$; see
\cite[Chaps.~7 and~9]{Kac90}. Write $A=(a_{ij})_{i,j\in\hI}$
for its Cartan matrix, where $\hI=\{0\}\sqcup \I$. Its simple roots are $\alpha_0=\delta-\theta$
and the finite simple roots $\alpha_j$, $j\in\I$, where $\theta$ is
the highest root of $\g$. We choose the positive simple generators so
that their images in the loop algebra are
\[
 \overline e_0=f_\theta v,\qquad
 \overline e_j=e_j\in\g_{\alpha_j}\quad(j\in\I),
\]
with $0\ne f_\theta\in\g_{-\theta}$.

\begin{lemma}\label{lem:identity-classical-relations}
With $v=u^2$, the above loop realization is an affine datum for
$(\g,\mathrm{id})$ with $e=2$, $s_0=2$ and $s_j=0$ for $j\in\I$. The
index set \eqref{eq:serre-index-set} is
$\mathcal S=\{(0,j):j\in\I,\ a_{0j}<0\}$, and the corresponding Serre
words are
\begin{equation}\label{eq:identity-classical-serre}
 r_j=(\ad E_0)^{1-a_{0j}}E_j,
 \qquad j\in\I,\quad a_{0j}<0,
\end{equation}
where $E_0=\underline{f_\theta}$ and $E_j=e_j$. Their degrees and
weights for $\mathfrak t=\mathfrak h$ are
\begin{equation}\label{eq:identity-relation-weights}
 d_j=2(1-a_{0j}),\qquad
 \mu_j=\alpha_j-(1-a_{0j})\theta,
\end{equation}
and $\mu_j\ne0$, $\g_{\mu_j}=0$ for every such $j$.
\end{lemma}

\begin{proof}
Since $\omega=\mathrm{id}$, the twisted loop algebra is
$L=\g[u^2,u^{-2}]=\g[v,v^{-1}]$, and (A1) holds with $s_0=2$ and
$s_j=0$. The elements $\overline h_i$ span $\mathfrak t=\mathfrak h$,
which generates $\kalg=\g$ together with the $e_j,f_j$, $j\in\I$. This
proves (A2). The lowest weight vector $f_\theta=y_0$ generates the
irreducible adjoint module under the finite positive simple
generators, which proves (A3).

The Serre words involving only finite nodes have degree zero. The
words $(\ad E_j)^{1-a_{j0}}E_0$, and $[E_0,E_j]$ when $a_{0j}=0$, have
degree $2=e$. Hence $\mathcal S$ consists of the pairs $(0,j)$ with
$a_{0j}<0$, and $r_j=r_{0j}$ has degree $2(1-a_{0j})>2$. The word $r_j$ contains $1-a_{0j}$ copies of
$E_0$, of weight $-\theta$, and one $E_j$, of weight $\alpha_j$, which
gives \eqref{eq:identity-relation-weights}. The last statement follows
from Lemma~\ref{lem:classical-relation-weights}, since $\g^{+}=\g$. It
can also be seen directly as follows. The element
$\beta_j=\alpha_j+(1-a_{0j})\alpha_0=(1-a_{0j})\delta+\mu_j$ is not an
affine root. On the other hand, in an untwisted affine root system,
$N\delta+\mu$ is a root whenever $\mu$ is a finite root, and
$N\delta$ is a root for $N>0$. Therefore, $\mu_j$ is neither a root nor
zero.
\end{proof}

For example, for $\g=\mathfrak{sl}_2$, we have $a_{01}=-2$ and
$\theta=\alpha_1$. Then $\mathcal S=\{(0,1)\}$, and the only Serre
word in the presentation is $(\ad\underline{f_1})^3e_1$, which has
degree six and weight $-2\alpha_1$. Since this weight does not occur
in the adjoint module, the argument also covers the rank one case.

\begin{theorem}\label{thm:identity-coideal-pbw}
For every $a\in\C$, $\B^{\mathrm{ev}}_a$ is
$\hbar$-saturated in $\Y$, and its natural reduction
is
\[
 \B^{\mathrm{ev}}_a/\hbar\B^{\mathrm{ev}}_a
       \cong\rU(\g[u^2])\subset\rU(\g[u]).
\]
Ordered monomials in homogeneous lifts of an ordered homogeneous basis
of $\g[u^2]$ form a $\C[\hbar]$-basis. In particular,
$\B^{\mathrm{ev}}_a$ is a homogeneous coideal quantization of
$\g[u^2]$.
\end{theorem}

\begin{proof}
We apply Proposition~\ref{prop:saturation-criterion} with $e=2$,
$V=\g$, the affine datum in Lemma~\ref{lem:identity-classical-relations}
and $\sfb=\mathcal G_a$. By \eqref{eq:identity-equivariance-limit},
the map $\mathcal G_a$ is $\g$-equivariant and satisfies (Q1).
Condition (Q2) is vacuous as $\omega=\mathrm{id}$, and (Q3) follows
from \eqref{eq:identity-generator-coproduct}, whose correction terms
have first factors in $\rU(\g)[\hbar]$.
\end{proof}

\subsection{Embedded classification and uniqueness}

Next, we compare an arbitrary coideal quantization of $\g[u^2]$ with
the algebras constructed above. Note that a degree-two lift is not
uniquely determined by the equivariance and the reduced coproduct,
but the embedded algebra determines a unique spectral parameter.

\begin{theorem}\label{thm:identity-uniqueness}
Every homogeneous coideal quantization $B\subset\Y$ of $\g[u^2]$
equals $\B^{\mathrm{ev}}_a$ for a unique $a\in\C$. For $a_0,a_1\in\C$,
the spectral translation restricts to a graded $\C[\hbar]$-algebra
isomorphism
\begin{equation}\label{eq:identity-parameter-isomorphism}
 \sigma_{a_1-a_0}:\B^{\mathrm{ev}}_{a_0}
       \xrightarrow{\sim}\B^{\mathrm{ev}}_{a_1},
 \qquad x\longmapsto x,\quad
 \mathcal G_{a_0}(x)\longmapsto\mathcal G_{a_1}(x).
\end{equation}
Thus the twisted Yangian associated with the identity involution
exists and is unique up to graded $\C[\hbar]$-algebra isomorphism.
\end{theorem}

\begin{proof}
By the classical identification, we have
\begin{equation}\label{eq:identity-low-degree-pieces}
 B_0=\rU(\g),\qquad B_1=\hbar\rU(\g).
\end{equation}
Indeed, the reduction is injective in degree zero, and the degree-one
component of $\rU(\g[u^2])$ vanishes, so $B_1\subset B\cap\hbar\Y_0
=\hbar B_0$ by saturation.

Let
\[
 \mathrm{cl}_2:B_2\longrightarrow\rU(\g[u^2])_2
\]
be the reduction modulo $\hbar$ followed by the classical
identification $\cli$, and choose a linear map $b:\g\to B_2$ such
that $\mathrm{cl}_2(b(x))=xu^2$ for all $x\in\g$. We do not require
$b$ to be $\g$-equivariant. Since $\g\subset B_0$ and $B$ is a
graded subalgebra, $B_2$ is stable under the adjoint action of $\g$.
Moreover, $\mathrm{cl}_2$ is $\g$-equivariant, as the reduction is
an algebra homomorphism and $\cli$ is the identity on $\g$. Hence, for
$z,x\in\g$, the element $[z,b(x)]-b([z,x])\in B_2$ has zero
reduction. By saturation and \eqref{eq:identity-low-degree-pieces}, we
have
\begin{equation}\label{eq:identity-lift-defect}
 [z,b(x)]-b([z,x])\in\hbar B_1=\hbar^2\rU(\g)
 \qquad(z,x\in\g).
\end{equation}

As $b(x)$ and $\mathcal G(x)$ have the same reduction, we may write
\[
 b(x)-\mathcal G(x)=\hbar S(x),
 \qquad S(x)\in\Y_1,
\]
where $S:\g\to\Y_1$ is linear. We now use the coideal property to
constrain $S$. By Lemma~\ref{lem:tensor-saturation}, $B\otimes\Y$ is
a saturated submodule of $\Y\otimes\Y$, and we regard it as such. Its
reduction modulo $\hbar$ is $\rU(\g[u^2])\otimes\rU(\g[u])$, and the
reduction of $\nabla(b(x))\in B\otimes\Y$ is the reduced coproduct of
$xu^2$, which vanishes because $xu^2$ is primitive. Hence
\[
 \nabla(b(x))=\hbar\,Q(x)
\]
for a unique $Q(x)\in B\otimes\Y$. Since $\nabla(b(x))$ is homogeneous
of total degree two, $Q(x)$ is homogeneous of total degree one, so
$Q(x)\in B_0\otimes\Y_1+B_1\otimes\Y_0$. By
\eqref{eq:identity-low-degree-pieces},
$B_1\otimes\Y_0=\hbar B_0\otimes\Y_0=B_0\otimes\hbar\Y_0\subseteq
B_0\otimes\Y_1$, and therefore
\[
 \hbar^{-1}\nabla(b(x))=Q(x)\in\rU(\g)\otimes\Y_1,
\]
where $\rU(\g)\otimes\Y_1$ denotes the image of
$\rU(\g)\otimes_{\C}\Y_1$ in $\Y\otimes\Y$. On the other hand, by
\eqref{eq:identity-generator-coproduct} with $a=0$, we have
\[
 \hbar^{-1}\nabla(\mathcal G(x))
 =-[1\otimes\mathcal J(x),\Omega_{\g}]+\hbar R(x)
 \in\rU(\g)\otimes\Y_1,
\]
since $[1\otimes\mathcal J(x),\Omega_{\g}]
=\sum_r x_r\otimes\mathcal J([x,x^r])$. Subtracting, we get
$\nabla(S(x))=\hbar^{-1}\bigl(\nabla(b(x))-\nabla(\mathcal G(x))\bigr)
\in\rU(\g)\otimes\Y_1$. Reducing modulo $\hbar$, we obtain
\[
 \Delta_0(\overline{S(x)})
   -\overline{S(x)}\otimes1-1\otimes\overline{S(x)}
       \in\rU(\g)\otimes\rU(\g[u]).
\]
Applying Lemma~\ref{lem:relative-primitives} to $\g\subset\g[u]$ with
complement $u\g[u]$, we obtain $\overline{S(x)}\in\rU(\g)+u\g[u]$.
As $\overline{S(x)}$ has degree one, we have
$\overline{S(x)}=A(x)u$ for a linear map $A:\g\to\g$.

We claim that $A$ is $\g$-equivariant, although $b$ need not be.
Indeed, since $\mathcal G$ is $\g$-equivariant by
Lemma~\ref{lem:identity-coideal-generators}, it follows from
\eqref{eq:identity-lift-defect} that
\[
 [z,S(x)]-S([z,x])
 =\hbar^{-1}\bigl([z,b(x)]-b([z,x])\bigr)\in\hbar\rU(\g)
 \qquad(z,x\in\g),
\]
where the division by $\hbar$ is unique because $\Y$ is torsion-free.
Reducing modulo $\hbar$, we obtain $[z,A(x)]u=A([z,x])u$, that is,
$[z,A(x)]=A([z,x])$. Since the adjoint module is irreducible, by
Schur's lemma, $A=2a\,\mathrm{id}_{\g}$ for some $a\in\C$.

The element $S(x)-2a\mathcal J(x)$ has degree one and zero reduction,
so it equals $\hbar T(x)$ with $T(x)\in\Y_0=\rU(\g)$. Thus
\begin{equation}\label{eq:identity-arbitrary-lift}
 b(x)=\mathcal G(x)+2a\hbar\mathcal J(x)
                    +\hbar^2T(x),
 \qquad T(x)\in\rU(\g).
\end{equation}
As $\rU(\g)\subset B$, we may subtract $\hbar^2T(x)$ and add
$a^2\hbar^2x$ inside $B$. Then we obtain $\mathcal G_a(x)\in B$ for
all $x\in\g$, and hence $\B^{\mathrm{ev}}_a\subseteq B$. By
Theorem~\ref{thm:identity-coideal-pbw}, the reduction of this
inclusion is the identity of $\rU(\g[u^2])$. Therefore,
$\B^{\mathrm{ev}}_a=B$ by Lemma~\ref{lem:graded-lifting}.

Now suppose that $\B^{\mathrm{ev}}_a=\B^{\mathrm{ev}}_{a'}$. In this
common subalgebra, we have
\[
 \mathcal G_a(x)-\mathcal G_{a'}(x)
 =2(a-a')\hbar\mathcal J(x)
        +(a^2-(a')^2)\hbar^2x,
\]
and the last term belongs to $\rU(\g)[\hbar]$. If $a\ne a'$, then by
saturation, $\mathcal J(x)\in B$ for every $x\in\g$. For $x\ne0$, this
contradicts $\overline{\mathcal J(x)}=xu$, since $\rU(\g[u^2])$ has no
degree-one component. Thus $a=a'$.

Finally, by \eqref{eq:identity-spectral-translation},
$\sigma_{a_1-a_0}$ maps the generators of $\B^{\mathrm{ev}}_{a_0}$ to
those of $\B^{\mathrm{ev}}_{a_1}$, and its inverse is
$\sigma_{a_0-a_1}$. This proves
\eqref{eq:identity-parameter-isomorphism}. The existence follows from
Theorem~\ref{thm:identity-coideal-pbw}.
\end{proof}

\begin{remark}
Note that there is no central parameter, since $\mathfrak z(\g)=0$.
The parameter $a$ is uniquely determined by the embedded subalgebra
and the fixed map $\mathcal J$, but the degree-two lift is not unique.
For example, for $\lambda\in\C$, $\mathcal G_a(x)+\lambda\hbar^2x$ is
another equivariant lift in the same subalgebra with the same reduced
coproduct, because $\nabla(\lambda\hbar^2x)=0$.
\end{remark}

\section{Finite presentations of split and quasi-split twisted Yangians}
\label{sec:minimal}

In this section, we give a presentation of the twisted Yangians
associated with split and quasi-split symmetric pairs, using a
degree-zero copy of $\kalg$ and a degree-one copy of $\palg$. Besides
the Lie algebra and module relations, we only need a family of
commutators between Cartan generators. We first specify these
relations and show that they hold in the quantizations of
Section~\ref{sec:parameter-families}. Then the completeness follows
from the classical presentation of the twisted current algebra and
Lemma~\ref{lem:graded-lifting}. For symmetric pairs that are not
quasi-split, we refer to \cite{BR17} for relations in the
$\mathcal J$-presentation.

\subsection{Lie algebra data and the defining relations}
\label{ssec:reference-parameter}

Let $\tau$ be an involution of the Dynkin diagram, extended to
$\g$ by permuting the Chevalley generators. Write
$\omega_{\mathrm{ch}}$ for the Chevalley involution,
\[
 \omega_{\mathrm{ch}}(e_i)=-f_i,\qquad
 \omega_{\mathrm{ch}}(f_i)=-e_i,\qquad
 \omega_{\mathrm{ch}}(h_i)=-h_i,
\]
and put
\begin{equation}\label{eq:split-quasi-split-involution}
 \omega=\omega_{\mathrm{ch}}\circ \tau.
\end{equation}
Then the pair $(\g,\g^{\omega})$ is of \emph{quasi-split type}, and it
is of \emph{split type} if $\tau=\mathrm{id}$. Throughout this
section, we assume either the split case in any finite type, or a
nontrivial quasi-split case in type $\mathsf A$, $\mathsf D$ or
$\mathsf E_6$, subject to
\begin{equation}\label{eq:presentation-scope}
 \g\not\cong\mathfrak{sl}_2\quad(\tau=\mathrm{id}),
 \qquad
 \g\not\cong\mathfrak{sl}_3\quad(\tau\ne\mathrm{id}).
\end{equation}
Note that these exclusions only concern the presentation in this
section; see \cite{Lu26min,HL26}. The existence and uniqueness
theorems hold for every symmetric pair.

Let $\mathfrak h$ and $\Phi^+$ be the Cartan subalgebra and the set of
positive roots. Choose dual root vectors
$x_\alpha^\pm\in\g_{\pm\alpha}$ with $(x_\alpha^+,x_\alpha^-)=1$ and
$\omega(x_\alpha^+)=-x_{\tau\alpha}^-$; for simple roots, we use the
normalization \eqref{eq:root-vector-normalization}. Put
\begin{equation}\label{eq:kb-lie-elements}
 \begin{aligned}
 \theta_i&=\xi_i-\xi_{\tau i}\in\kalg,&
 t_i&=\xi_i+\xi_{\tau i}\in\palg,\\
 b_\alpha&=x_\alpha^+-x_{\tau\alpha}^-\in\kalg,&
 y_\alpha&=x_\alpha^++x_{\tau\alpha}^-\in\palg.
 \end{aligned}
\end{equation}
The $\theta_i$ span $\mathfrak h\cap\kalg$, and the $t_i$ span
$\mathfrak h\cap\palg$. Together with the $b_\alpha$, respectively the
$y_\alpha$, they span $\kalg$, respectively $\palg$. Set
\[
 a_{i\alpha}:=(\alpha,\alpha_i+\alpha_{\tau i}),\qquad \{v,w\}=vw+wv.
\]
Then $[t_i,b_\alpha]=a_{i\alpha}y_\alpha$ and
$a_{i,\tau\alpha}=a_{i\alpha}$.

We take $\K_\circ=0$, except in nontrivial quasi-split type
$\mathsf A_{2n}$. In that type, the Cartan relations below hold at a
nonzero central parameter $\K_\circ$ and require a linear correction
$\eta_i$. We now specify both of them, following \cite{HL26}. Put
$N=2n+1$, so that
$\g=\mathfrak{sl}_N$ and $\tau i=N-i$. Let $E_{ij}$ be the matrix
units, put $i'=N+1-i$, and use the form $(x,y)=\operatorname{tr}(xy)$.
For $i<j$, choose
\begin{equation}\label{eq:A2n-root-vectors}
 x_{\varepsilon_i-\varepsilon_j}^+
   =(\sqrt{-1})^{j-i-1}E_{ij},\qquad
 x_{\varepsilon_i-\varepsilon_j}^-
   =(\sqrt{-1})^{i-j+1}E_{ji},
\end{equation}
where $\varepsilon_i$ is the $i$th diagonal coordinate functional.
The two involutions are
\begin{equation}\label{eq:A2n-involutions}
 \tau(E_{ij})=(-1)^{i-j+1}E_{j'i'},\qquad
 \omega(E_{ij})=(-1)^{i-j}E_{i'j'}.
\end{equation}
Define
\begin{equation}\label{eq:A2n-central-element}
 \mathscr S=\sum_{i=1}^{N}(-1)^{n-i}E_{i'i},\qquad
 \mathscr Z=\mathscr S+N^{-1}1_N,\qquad
 \K_\circ=\tfrac18\mathscr Z\in\mathfrak z(\kalg),
\end{equation}
see \cite[Lem.~2.9]{HL26}. Thus $\mathscr Z$ is the element
$\mathcal K$ of \cite[(2.14)]{HL26}, and $\B_{0,\K_\circ}$ is the
twisted Yangian of \cite[Def.~2.8]{HL26} with $c=\tfrac18$.
For this type, we set
\[
 \wp_i=\begin{cases}
  1,&i=n,\\
  -1,&i=n+1,\\
  0,&\text{otherwise},
 \end{cases}
 \qquad \kappa_i=\varepsilon_i-\varepsilon_{i'},
 \qquad
 \eta_i=\tfrac14\wp_i\theta_i+
          \tfrac18\sqrt{-1}\sum_{j=1}^n a_{i\kappa_j}b_{\kappa_j}.
\]
In all other types, we set $\wp_i=0$ and $\eta_i=0$. Hence
$\eta_i\in\kalg$ in all cases. Define
\begin{equation}\label{eq:kb-q}
 q_i=\tfrac12\theta_i^2+
       \tfrac14\sum_{\alpha\in\Phi^+}(\alpha,\alpha_i)
                    \{b_\alpha,b_{\tau\alpha}\}+\eta_i
       \in\rU(\kalg).
\end{equation}
Note that the $q_i$ have PBW degree at most two and current degree
zero. In the split case, they simplify to
$q_i=\tfrac12\sum_{\alpha\in\Phi^+}(\alpha,\alpha_i)b_\alpha^2$.

\begin{remark}\label{rem:q-tau-invariance}
Since $t_{\tau i}=t_i$, the relations \eqref{eq:kb-extra} below for the
indices $i$ and $\tau i$ must be compatible. Indeed, we have
$q_{\tau i}=q_i$. To see this, note that $\theta_{\tau i}=-\theta_i$.
The substitution $\alpha\mapsto\tau\alpha$ identifies the root parts
of $q_i$ and $q_{\tau i}$, because
$(\tau\alpha,\alpha_{\tau i})=(\alpha,\alpha_i)$ and
$\{b_{\tau\alpha},b_\alpha\}=\{b_\alpha,b_{\tau\alpha}\}$. Moreover,
$\wp_{\tau i}\theta_{\tau i}=\wp_i\theta_i$ and
$a_{\tau i,\kappa}=a_{i\kappa}$.
\end{remark}

Let $\mathcal U^\imath_\hbar(\g,\tau)$ be the unital associative
$\C[\hbar]$-algebra generated by the images of two linear maps
\[
 \iota:\kalg\longrightarrow\mathcal U^\imath_\hbar(\g,\tau)_0,
 \qquad
 \mathcal B:\palg\longrightarrow\mathcal U^\imath_\hbar(\g,\tau)_1,
 \qquad \deg\hbar=1,
\]
subject to the following relations:
\begin{enumerate}[label=(\arabic*)]
\item $\iota$ is a Lie algebra homomorphism and $\mathcal B$ is a
$\kalg$-module homomorphism:
\begin{equation}\label{eq:kb-equivariance}
 [\iota(x),\iota(x')]=\iota([x,x']),\qquad
 [\iota(x),\mathcal B(y)]=\mathcal B([x,y])
 \quad(x,x'\in\kalg,\ y\in\palg).
\end{equation}
We also write $\iota:\rU(\kalg)\to\mathcal U^\imath_\hbar(\g,\tau)$
for its unital extension.
\item For all $i,j\in\I$,
\begin{equation}\label{eq:kb-extra}
 [\mathcal B(t_i),\mathcal B(t_j)]
       =\hbar^2\iota([q_j,q_i]).
\end{equation}
\end{enumerate}
We call \eqref{eq:kb-extra} the \emph{Cartan relations}. After choosing
bases of $\kalg$ and $\palg$, we obtain a finite list of generators
and relations. Here the elements $q_i$ are fixed polynomials in
the degree-zero generators, and they are not additional generators.

\begin{remark}
Under \eqref{eq:presentation-scope}, the extra low-rank relations in
\cite[Thm.~3.1]{Lu26min} and \cite[Thm.~3.1 and Cor.~4.13]{HL26} do not occur, and the above presentation
is a reformulation of the minimalistic presentation there. Thus
Theorem~\ref{thm:main-presentation} below gives a direct proof of
these presentations in the $\mathcal J$-presentation, starting from
the embedding into the Yangian.
\end{remark}

\subsection{Cartan relations}
By Lemma~\ref{lem:lift}, the relations \eqref{eq:kb-equivariance} hold
for the normalized lift of every quantization. To prove
\eqref{eq:kb-extra} at $\K_\circ$, we show that the corrected elements
$\mathcal B(t_i)+\hbar\iota(q_i)$ commute. We first show that the
mixed terms in their commutator cancel.

\begin{lemma}\label{lem:kb-mixed-terms}
Let $A$ be a $\C[\hbar]$-algebra, $\iota:\kalg\to A$ a Lie algebra
homomorphism, and $\mathcal B:\palg\to A$ a $\kalg$-module
homomorphism. Then
\begin{equation}\label{eq:kb-cross-cancellation}
 [\mathcal B(t_i),\iota(q_j)]
       =[\mathcal B(t_j),\iota(q_i)]\qquad(i,j\in\I).
\end{equation}
\end{lemma}

\begin{proof}
By equivariance, we have
\[
 [\mathcal B(t_i),\iota(\theta_j)]=0,\qquad
 [\mathcal B(t_i),\iota(b_\alpha)]
                  =a_{i\alpha}\mathcal B(y_\alpha).
\]
When computing $[\mathcal B(t_i),\iota(q_j)]$, we reindex the terms
involving $[\mathcal B(t_i),\iota(b_{\tau\alpha})]$ by
$\alpha\mapsto\tau\alpha$. As $a_{i,\tau\alpha}=a_{i\alpha}$ and
$(\alpha,\alpha_j)+(\tau\alpha,\alpha_j)=a_{j\alpha}$, we obtain
\[
 [\mathcal B(t_i),\iota(q_j)]
 =\tfrac14\sum_{\alpha>0}a_{i\alpha}a_{j\alpha}
       \{\mathcal B(y_\alpha),\iota(b_{\tau\alpha})\}
       +\mathcal B([t_i,\eta_j]).
\]
Clearly the sum is symmetric in $i$ and $j$. The last term vanishes
outside nontrivial quasi-split type $\mathsf A_{2n}$, and in that type,
we have
\[
 [t_i,\eta_j]=\tfrac18\sqrt{-1}\sum_{l=1}^n
                  a_{i\kappa_l}a_{j\kappa_l}y_{\kappa_l},
\]
which is also symmetric.
\end{proof}

The following observation will be used to show that the corrected
commutator vanishes. Here weights are taken with respect to the
Cartan subalgebra
$\mathfrak h\subset\g$. Put $Q_+=\sum_{i\in\I}\mathbb Z_{\geqslant0}\alpha_i$
and
\[
 \Y_{>0}=\bigoplus_{\beta\in Q_+\setminus\{0\}}\Y_\beta,
 \qquad
 \Y_\beta=\{v\in\Y:[h,v]=\beta(h)v\text{ for all }h\in\mathfrak h\}.
\]

\begin{lemma}\label{lem:positive-weight-separation}
For every homogeneous coideal quantization $B\subset\Y$ of the pair
under consideration,
\begin{equation}\label{eq:positive-weight-separation}
 B\cap\Y_{>0}=0.
\end{equation}
\end{lemma}

\begin{proof}
Write $\mathfrak b^-=\mathfrak n^-\oplus\mathfrak h$ and order a basis
of $\g[u]$ with $\mathfrak b^-[u]$ before $\mathfrak n^+[u]$. By the
PBW theorem, we have a linear projection
\[
 \varpi:\rU(\g[u])\longrightarrow\rU(\mathfrak b^-[u])
\]
which keeps the ordered monomials with no factor in
$\mathfrak n^+[u]$. We first show that its restriction to $\rU(\Ltw)$
is injective. The corresponding projection
$\varpi_1:\g[u]\to\mathfrak b^-[u]$ is injective on $\Ltw$. Indeed, an
element of its kernel lies in $\Ltw\cap\mathfrak n^+[u]$. It is fixed
by $\check\omega$, and $\check\omega$ maps it into $\mathfrak n^-[u]$,
so it is zero. For the PBW filtrations, the top-degree part of
$\varpi$ on $\rU(\Ltw)$ is the injective map
$\rS(\varpi_1|_{\Ltw}):\rS(\Ltw)\to\rS(\mathfrak b^-[u])$. Hence the
highest PBW term of a nonzero element of $\rU(\Ltw)$ survives, and
$\varpi|_{\rU(\Ltw)}$ is injective.

The projection $\varpi$ maps each $\mathfrak h$-weight space into the
weight space of the same weight. As its target has only weights in
$-Q_+$, it vanishes on every strictly positive weight space.
Consequently,
\[
 \rU(\Ltw)\cap
 \bigoplus_{\beta\in Q_+\setminus\{0\}}\rU(\g[u])_\beta=0.
\]
Let $v\in B\cap\Y_{>0}$ be of degree $d$. Then its reduction lies in
this zero intersection, and by saturation, $v=\hbar v_1$ for some
$v_1\in B_{d-1}$. Since $\hbar$ has weight zero and multiplication by
$\hbar$ is injective, we have $v_1\in\Y_{>0}$. Repeating this argument,
we lower the degree and obtain $v=0$. The case of arbitrary elements
follows by taking homogeneous components.
\end{proof}

The lemma concerns elements all of whose weight components
are strictly positive. Individual elements of $B$ may have positive
weight components, and $B$ need not be stable under $\mathfrak h$.

\begin{lemma}\label{lem:kb-cancellation}
The normalized lift $\mathcal B$ of every homogeneous coideal
quantization with parameters $(a,\K_\circ)$, together with the
inclusion $\iota:\kalg\to B$, satisfies the Cartan relations
\eqref{eq:kb-extra}.
\end{lemma}

\begin{proof}
Put $\mathfrak Z_i=\mathcal B(t_i)+\hbar q_i$. By
Lemma~\ref{lem:kb-mixed-terms},
\begin{equation}\label{eq:kb-cartan-defect}
 [\mathfrak Z_i,\mathfrak Z_j]=[\mathcal B(t_i),\mathcal B(t_j)]+\hbar^2[q_i,q_j],
\end{equation}
so it remains to show that $[\mathfrak Z_i,\mathfrak Z_j]=0$. For
this, we separate each $\mathfrak Z_i$ into a weight-zero part and a
strictly positive-weight part.

For $h\in\mathfrak h$, set
\[
 \mathfrak W(h)=\tfrac14\sum_{\alpha>0}\alpha(h)
                         \{x_\alpha^+,x_\alpha^-\},\qquad
 \mathcal H(h)=\mathcal J(h)-\hbar \mathfrak W(h).
\]
We claim that
\begin{equation}\label{eq:ambient-commuting-cartan-polynomials}
 [\mathcal H(h),\mathcal H(h')]=0,\qquad [h,\mathcal H(h')]=0
 \qquad(h,h'\in\mathfrak h);
\end{equation}
see e.g.\ \cite[Thm.~2.6]{GRW19}. Indeed, recall $\nu(h)$ from
\eqref{eq:nu}. Note that the difference
$$\mathfrak W(h)-\nu(h)
=\tfrac14\sum_{\alpha>0}\alpha(h)[x_\alpha^+,x_\alpha^-]$$ lies in
$\mathfrak h$, so it commutes with $\mathcal J(\mathfrak h)$ and with
the weight-zero elements $\nu(\mathfrak h)$. Moreover,
\[
 [\mathcal J(h),\nu(h')]=\tfrac12\sum_{\alpha>0}\alpha(h)\alpha(h')
 \bigl(x_\alpha^-\mathcal J(x_\alpha^+)-\mathcal J(x_\alpha^-)x_\alpha^+\bigr)
\]
is symmetric in $h,h'$. Hence the first identity follows from the
defining relation
\[
 [\mathcal J(\xi_i),\mathcal J(\xi_j)]=\hbar^2[\nu(\xi_j),\nu(\xi_i)],
\]
and the second one holds because $\mathcal J(h')$ and
$\mathfrak W(h')$ have weight zero.

The root-vector part of $\kalg$ is orthogonal to $\mathfrak h\cap\kalg$,
and $(b_\alpha,b_\beta)=-2\delta_{\beta,\tau\alpha}$. Thus, writing
$C_{\mathfrak h\cap\kalg}$ for the Cartan part of the Casimir element,
we have
\[
 C_{\kalg}=C_{\mathfrak h\cap\kalg}
                 -\tfrac12\sum_{\alpha>0}b_\alpha b_{\tau\alpha}.
\]
As $t_i$ commutes with the Cartan part, expanding the bracket and
reindexing by $\tau$, we obtain
\[
 -\tfrac14[t_i,C_{\kalg}]
       =\tfrac18\sum_{\alpha>0}a_{i\alpha}
                                  \{y_\alpha,b_{\tau\alpha}\}.
\]
Similarly, the root part of $q_i$ is equal to
$\tfrac18\sum_{\alpha>0}a_{i\alpha}\{b_\alpha,b_{\tau\alpha}\}$.
Adding these and using $y_\alpha+b_\alpha=2x_\alpha^+$ and
$a_{i\alpha}=\alpha(t_i)$, we obtain
\begin{equation}\label{eq:cartan-triangular-expansion}
 -\tfrac14[t_i,C_{\kalg}]+q_i
     =\tfrac12\theta_i^2-\mathfrak W(t_i)+\mathfrak P_i+\eta_i,
 \qquad
 \mathfrak P_i=\tfrac14\sum_{\alpha>0}a_{i\alpha}
                            \{x_\alpha^+,x_{\tau\alpha}^+\}.
\end{equation}
Each term of $\mathfrak P_i$ has weight
$\alpha+\tau\alpha\in Q_+\setminus\{0\}$. Therefore, by the embedding
\eqref{eq:embedding-rigid} of the normalized lift with
$z=z_{a,\K_\circ}$, we have
\begin{equation}\label{eq:normalized-cartan-triangular-form}
 \mathfrak Z_i=\mathcal H(t_i)+\hbar\bigl(at_i+\tfrac12\theta_i^2\bigr)
            +\hbar \mathfrak P_i+\hbar\bigl(\eta_i+[\K_\circ,t_i]\bigr).
\end{equation}
The last term in parentheses vanishes, except in nontrivial
quasi-split type $\mathsf A_{2n}$.

In type $\mathsf A_{2n}$, we put $\mathfrak R_i=(-1)^{n-i}E_{ii'}$. By
the matrix formulas above, we have
\[
 \omega(\mathfrak R_i)=(-1)^{n-i}E_{i'i},\qquad
 b_{\kappa_i}=\sqrt{-1}\bigl(\mathfrak  R_i+\omega(\mathfrak R_i)\bigr).
\]
The off-diagonal part of $\mathscr Z$ is
$\sum_{1\leqslant i\leqslant n}(\mathfrak R_i+\omega(\mathfrak R_i))$,
and its diagonal part commutes with $t_i$. Hence
\[
 [\mathscr Z,t_i]=-
        \sum_{j=1}^n a_{i\kappa_j}\bigl(\mathfrak R_j-\omega(\mathfrak R_j)\bigr),
 \qquad
 \eta_i=\tfrac14\wp_i\theta_i-
        \tfrac18\sum_{j=1}^n a_{i\kappa_j}\bigl(\mathfrak R_j+\omega(\mathfrak R_j)\bigr).
\]
Since $\K_\circ=\tfrac18\mathscr Z$, we obtain
\[
 \eta_i+[\K_\circ,t_i]
       =\tfrac14\wp_i\theta_i-
          \tfrac14\sum_{j=1}^n a_{i\kappa_j}\mathfrak R_j.
\]
Thus the negative-root terms cancel, and each $\mathfrak R_j$ has
positive weight $\kappa_j$.

In all cases, we have shown that
\[
 \mathfrak Z_i=\mathfrak Z_i^0+\hbar\mathfrak Z_i^{+},\qquad
 \mathfrak Z_i^0=\mathcal H(t_i)+\hbar
       \bigl(at_i+\tfrac12\theta_i^2+\tfrac14\wp_i\theta_i\bigr),
 \qquad \mathfrak Z_i^{+}\in\Y_{>0}.
\]
The elements $\mathfrak Z_i^0$ have weight zero, and by
\eqref{eq:ambient-commuting-cartan-polynomials}, they commute with
each other, as $t_i$, $\theta_i^2$ and $\theta_i$ lie in
$\rU(\mathfrak h)$. In the expansion of $[\mathfrak Z_i,\mathfrak Z_j]$,
every remaining term has at least one factor in $\Y_{>0}$, and all
other factors have weight zero or lie in $\Y_{>0}$. Hence every
remaining term lies in $\Y_{>0}$. Thus
$[\mathfrak Z_i,\mathfrak Z_j]\in B\cap\Y_{>0}$, which is zero by
Lemma~\ref{lem:positive-weight-separation}. Now the lemma follows from
\eqref{eq:kb-cartan-defect}.
\end{proof}

\subsection{Completeness of the presentation}
\label{ssec:arbitrary-central-presentation}

We use the classical presentation of the twisted current algebra in
the language of \S\ref{ssec:free-model}, with $e=1$ and $V=\palg$.
Under \eqref{eq:presentation-scope}, the kernel of
$\rU(\mathfrak F)\to\rU(\Ltw)$, $x\mapsto x$,
$\underline y\mapsto yu$, is generated as a two-sided ideal by
\begin{equation}\label{eq:classical-k-p-cartan-relation}
 [\underline{t_i},\underline{t_j}]\qquad(i,j\in\I).
\end{equation}
This is the minimalistic presentation in \cite[Prop.~3.2]{Lu26min} for
the split case, where the redundant relations in types $\mathsf B_2$
and $\mathsf G_2$ are removed by \cite[Lem.~4.12]{HL26}, and the one
in \cite[Prop.~3.2]{HL26} for the nontrivial quasi-split case. We only
use these classical results here.

Note that although \eqref{eq:classical-k-p-cartan-relation} only
involves Cartan elements, the degree-one generating space is the full
module $\palg$, since the $t_i$ generate $\palg$ as a $\kalg$-module.
Indeed, choose $h\in\mathfrak h$ such that $\alpha_i(h)>0$ for every
simple root, and put $t=\tfrac12(h+\tau(h))$. Then $\tau(t)=t$, so
$t\in\operatorname{span}_{\C}\{t_i:i\in\I\}$. As $\tau$ preserves
positive roots, we have
\[
 \alpha(t)=\tfrac12\bigl(\alpha(h)+(\tau\alpha)(h)\bigr)>0,
 \qquad
 [t,b_\alpha]=\alpha(t)y_\alpha
 \qquad(\alpha\in\Phi^+).
\]
Thus the $\kalg$-submodule generated by the $t_i$ contains every
$y_\alpha$, and the $y_\alpha$ together with the $t_i$ span $\palg$.

\begin{theorem}\label{thm:main-presentation}
Under \eqref{eq:presentation-scope}, for every $a\in\C$ the
assignments
\begin{equation}\label{eq:central-parameter-embedding}
 \iota(x)\longmapsto x,\qquad
 \mathcal B(y)\longmapsto
 \mathcal J(y)-\tfrac14\hbar[y,C_{\kalg}]
                  +\hbar\bigl(ay+[\K_\circ,y]\bigr)
\end{equation}
define a graded algebra isomorphism
\[
 \mathcal U^\imath_\hbar(\g,\tau)
             \xrightarrow{\sim}\B_{a,\K_\circ}.
\]
In particular, $\mathcal U^\imath_\hbar(\g,\tau)$ is $\C[\hbar]$-free
with classical limit
\begin{equation}\label{eq:central-classical-marking}
 \mathcal U^\imath_\hbar(\g,\tau)/\hbar
             \mathcal U^\imath_\hbar(\g,\tau)
       \xrightarrow{\sim}\rU(\Ltw),\qquad
 \overline{\iota(x)}\longmapsto x,\quad
 \overline{\mathcal B(y)}\longmapsto yu,
\end{equation}
and for any ordered homogeneous basis of $\Ltw$, the ordered monomials
in homogeneous lifts of its elements form a $\C[\hbar]$-basis of the
presented algebra.
\end{theorem}

\begin{proof}
Put $A=\mathcal U^\imath_\hbar(\g,\tau)$ and $B=\B_{a,\K_\circ}$. By
Lemmas~\ref{lem:lift} and~\ref{lem:kb-cancellation}, the normalized
generators of $B$ satisfy the relations \eqref{eq:kb-equivariance}
and the Cartan relations \eqref{eq:kb-extra}, respectively. Therefore,
evaluating the two generating maps, we obtain a graded surjective
homomorphism $f:A\to B$.

We first show that $\overline f:A/\hbar A\to B/\hbar B$ is an
isomorphism. Modulo $\hbar$, the relations
\eqref{eq:kb-equivariance}--\eqref{eq:kb-extra} reduce to the Lie
algebra and module relations together with
\eqref{eq:classical-k-p-cartan-relation}. Hence the classical
presentation gives a surjective homomorphism
\[
 \psi:\rU(\Ltw)\longrightarrow A/\hbar A,\qquad
 x\longmapsto\overline{\iota(x)},\quad
 yu\longmapsto\overline{\mathcal B(y)},
\]
since these images generate $A/\hbar A$. On the other hand, by
Theorem~\ref{thm:coideal-pbw}, we have the isomorphism
$\cli:B/\hbar B\to\rU(\Ltw)$, and the map
$\Xi=\cli\circ\overline f:A/\hbar A\to\rU(\Ltw)$ satisfies
$\Xi(\overline{\iota(x)})=x$ and $\Xi(\overline{\mathcal B(y)})=yu$.
Thus $\Xi\circ\psi$ fixes the generators $\kalg$ and $u\palg$ of
$\rU(\Ltw)$, so $\Xi\circ\psi=\mathrm{id}$. In particular, $\psi$ is
injective. Therefore, $\psi$ is an isomorphism with inverse $\Xi$, and
$\overline f=\cli^{-1}\circ\psi^{-1}$ is an isomorphism.

As $B$ is a subalgebra of $\Y$, it is torsion-free. Hence $f$ is an
isomorphism by Lemma~\ref{lem:graded-lifting}. We do not need any
freeness assumption on $A$. The statements about the
freeness and the PBW basis follow from
Theorem~\ref{thm:coideal-pbw}.
\end{proof}

By Theorem~\ref{thm:embedded-classification}, this presentation applies
to every homogeneous coideal quantization with the corresponding
parameters. Presentations for other central parameters can be
obtained by applying the central translation in
Lemma~\ref{lem:central-translation-quantization}.

\appendix
\section{The technical assumption for quantum duality}
\label{sec:quantum-duality-assumption}

In this appendix, we verify the technical assumptions
\cite[(3.5) and (3.13)]{LWW25} for the Yangian and for every
homogeneous coideal quantization considered in this paper. Throughout
this appendix, $\omega$ is an arbitrary involution, including the
identity, and tensor products are taken over $\C[\hbar]$. We use the
right coideal convention, and explain how to pass to the left
convention at the end.

Recall the counit of $\Y$ from \eqref{eq:couint}. Let
$\Delta^{(n)}:\Y\to\Y^{\otimes n}$ be the iterated coproduct, with
$\Delta^{(1)}=\mathrm{id}$, and put
\[
 \delta_n=(\mathrm{id}-\varepsilon)^{\otimes n}
                    \circ\Delta^{(n)}
 \qquad(n\geqslant1).
\]
The subalgebra used in the quantum duality principle
\cite{G02,CG06} is
\begin{equation}\label{eq:duality-ambient-prime}
 \Y'=\{X\in\Y:\delta_n(X)\in
             \hbar^n\Y^{\otimes n}\text{ for all }n\geqslant1\}.
\end{equation}
For a homogeneous right coideal quantization $B\subset\Y$, the
corresponding subalgebra is
\begin{equation}\label{eq:duality-coideal-prime}
 B'=\{b\in B:\delta_n(b)\in
       \hbar^n(B\otimes\Y^{\otimes(n-1)})
                    \text{ for all }n\geqslant1\}.
\end{equation}
Note that the condition for $n=0$ in \cite{LWW25} holds
automatically. The assumption \cite[(3.13)]{LWW25} requires
homogeneous lifts $b_\lambda\in B$ of a basis of $\Ltw$ such that
$\hbar b_\lambda\in B'$, and the analogous condition for $\g[u]$ and
$\Y'$ is \cite[(3.5)]{LWW25}.

\begin{lemma}\label{lem:duality-prime-properties}
The subspace $\Y'$ is a $\C[\hbar]$-subalgebra and satisfies
$[\Y',\Y']\subset\hbar\Y'$. For every homogeneous coideal quantization
$B\subset\Y$,
\begin{equation}\label{eq:duality-prime-intersection}
 B'=B\cap\Y'.
\end{equation}
Moreover, the $\C[\hbar]$-submodule
\[
 \mathcal A=\{X\in\Y:\hbar X\in\Y'\}
\]
is closed under Lie brackets, and $\hbar XY\in\mathcal A$ whenever
$X,Y\in\mathcal A$.
\end{lemma}

\begin{proof}
For $I\subseteq\{1,\ldots,n\}$, let $\delta_I(X)$ denote
$\delta_{|I|}(X)$ placed in the positions indexed by $I$, with units
in the other positions. For $I=\varnothing$, we put
$\delta_{\varnothing}(X)=\varepsilon(X)1^{\otimes n}$. By the counit
identities, we have $\Delta^{(n)}(X)=\sum_{I}\delta_I(X)$. Since
$\varepsilon$ is an algebra homomorphism, applying the reduced
projection to a product, we obtain
\[
 \delta_n(XY)
 =\sum_{I\cup J=\{1,\ldots,n\}}
                    \delta_I(X)\delta_J(Y).
\]
For $X,Y\in\Y'$, each summand is divisible by $\hbar^{|I|+|J|}$ and
$|I|+|J|\geqslant n$, so $XY\in\Y'$. Similarly,
\[
 \delta_n([X,Y])
 =\sum_{I\cup J=\{1,\ldots,n\}}
                    [\delta_I(X),\delta_J(Y)].
\]
The terms with $I\cap J=\varnothing$ vanish, because their tensor
positions are disjoint (this includes the case $I=\varnothing$ or
$J=\varnothing$, where one factor is a scalar). Every remaining term is
divisible by $\hbar^{|I|+|J|}$ with $|I|+|J|\geqslant n+1$. Moreover,
$X-\varepsilon(X)1$ and $Y-\varepsilon(Y)1$ lie in $\hbar\Y$, so
$[X,Y]/\hbar\in\Y$. By the divisibility shown above, we obtain
$[X,Y]/\hbar\in\Y'$.

To prove \eqref{eq:duality-prime-intersection}, note that by the
coideal property and the counit, we have
$\delta_n(B)\subset B\otimes\Y^{\otimes(n-1)}$. Since $B$ is saturated
in $\Y$ and $\Y^{\otimes(n-1)}$ is free, by
Lemma~\ref{lem:tensor-saturation}, $B\otimes\Y^{\otimes(n-1)}$ is
saturated in $\Y^{\otimes n}$. Hence by repeated division, we have
\[
 (B\otimes\Y^{\otimes(n-1)})\cap
       \hbar^n\Y^{\otimes n}
 =\hbar^n(B\otimes\Y^{\otimes(n-1)}).
\]
Finally, if $X,Y\in\mathcal A$, then
\[
 \hbar[X,Y]=\frac{[\hbar X,\hbar Y]}{\hbar}\in\Y',
 \qquad
 \hbar(\hbar XY)=(\hbar X)(\hbar Y)\in\Y'.\qedhere
\]
\end{proof}

The statement for Yangians in the following theorem is also recorded
in \cite[Rem.~3.2]{LWW25}. We include a direct proof together with the
coideal statement.

\begin{theorem}\label{thm:duality-technical-assumption}
Every homogeneous basis of $\g[u]$ admits homogeneous lifts
$X_\gamma\in\Y$ with $\hbar X_\gamma\in\Y'$. For every homogeneous
coideal quantization $B\subset\Y$ of $\Ltw$, every homogeneous basis
of $\Ltw$ admits homogeneous lifts $b_\lambda\in B$ with
$\hbar b_\lambda\in B'$. Thus the assumptions
\cite[(3.5) and (3.13)]{LWW25} hold for every involution of a complex
simple Lie algebra and every embedded quantization in the fixed
Yangian $\Y$.
\end{theorem}

\begin{proof}
We first consider the generators of the Yangian. For $x\in\g$, since
$x$ is primitive, we have $\delta_n(\hbar x)=0$ for $n\geqslant2$. By
the coproduct \eqref{eq:coproduct-J}, we have
\[
 \delta_1(\hbar\mathcal J(x))=\hbar\mathcal J(x),\qquad
 \delta_2(\hbar\mathcal J(x))
       =\tfrac12\hbar^2[x\otimes1,\Omega_{\g}],\qquad
 \delta_n(\hbar\mathcal J(x))=0\quad(n\geqslant3);
\]
for the last equality, note that every term of the iterated coproduct
of $\mathcal J(x)$ has at most two nonunit tensor factors, because the
correction term lies in $\hbar(\g\otimes\g)$ and $\g$ is primitive.
Hence $x,\mathcal J(x)\in\mathcal A$. Since $\g[u]$ is generated as a
Lie algebra by $\g$ and $u\g$, every homogeneous element of $\g[u]$ is
a linear combination of homogeneous Lie words in these generators.
Replacing $xu$ by $\mathcal J(x)$, we obtain a homogeneous lift, which
lies in $\mathcal A$ by Lemma~\ref{lem:duality-prime-properties}.

Now suppose that $\omega\ne\mathrm{id}$. By
Theorem~\ref{thm:embedded-classification}, $B=\B_{a,\K}$ for some $a$
and $\K$. For dual bases $v_r,v^r$ of $\kalg$, we have
\[
 [y,C_{\kalg}]
 =\sum_r\bigl([y,v_r]v^r+v_r[y,v^r]\bigr)
 \qquad(y\in\palg),
\]
and every factor on the right-hand side lies in $\g\subset\mathcal A$.
Thus $\hbar[y,C_{\kalg}]\in\mathcal A$ by the product property in
Lemma~\ref{lem:duality-prime-properties}. The remaining terms of
$b_{a,\K}(y)$ also lie in $\mathcal A$, so $\kalg$ and
$b_{a,\K}(\palg)$ lie in $B\cap\mathcal A$. By
Lemma~\ref{lem:bracket-generation}, $\kalg$ and $u\palg$ generate
$\Ltw$. We express the given basis vectors of $\Ltw$ as homogeneous
Lie words in these generators, and replace the generators by the above
lifts. Then the resulting elements $b_\lambda$ lie in
$B\cap\mathcal A$, and hence $\hbar b_\lambda\in B\cap\Y'=B'$.

For $\omega=\mathrm{id}$, we have $B=\B^{\mathrm{ev}}_a$ by
Theorem~\ref{thm:identity-uniqueness}. The element $\mathcal Q(x)$ in
\eqref{eq:identity-degree-two-generator} is a linear combination of
brackets of elements of $\mathcal J(\g)\subset\mathcal A$, so
$\mathcal Q(x)\in\mathcal A$. Moreover,
\[
 [\mathcal J(x),C_{\g}]
 =\sum_r\bigl(\mathcal J([x,x_r])x^r
                      +x_r\mathcal J([x,x^r])\bigr),
\]
so $\hbar[\mathcal J(x),C_{\g}]\in\mathcal A$, and therefore
$\mathcal G_a(x)\in\mathcal A$. As $\g[u^2]$ is generated by $\g$
and $u^2\g$, the same construction with Lie words gives the required
lifts.
\end{proof}

Note that the proof works for all spectral and central parameters, and
it also gives the left coideal statement used in \cite{LWW25}. Indeed,
the Yangian anti-involution $\mathfrak T$ in \S\ref{ssec:left-coideal}
preserves the counit and satisfies
\[
 \delta_n(\mathfrak T X)
 =\operatorname{rev}_n\,
       \mathfrak T^{\otimes n}\delta_n(X),
\]
where $\operatorname{rev}_n$ reverses the tensor factors. Thus
$\mathfrak T(\Y')=\Y'$, and $\mathfrak T$ sends the above lifts to
lifts satisfying the left coideal version of
\eqref{eq:duality-coideal-prime}.

\bibliographystyle{amsalpha}
\bibliography{reference}

\end{document}